\documentclass{amsart}
\usepackage{graphicx}
\usepackage{amscd}
\usepackage{amsmath}
\usepackage{amsxtra}
\usepackage{amsfonts}
\usepackage{amssymb}

\newtheorem{theorem}{Theorem}[section]
\newtheorem{corollary}[theorem]{Corollary}
\newtheorem{lemma}[theorem]{Lemma}
\newtheorem{proposition}[theorem]{Proposition}

\newtheorem{conjecture}[theorem]{Conjecture}
\theoremstyle{definition}
\newtheorem{definition}[theorem]{Definition}
\newtheorem{remark}[theorem]{Remark}
\newtheorem{problem}{Problem}

\newtheorem{example}[theorem]{Example}
\theoremstyle{remark}

\renewcommand{\theclaim}{\textup{\theclaim}}

\newtheorem*{acknowledgements}{Acknowledgements}

\numberwithin{equation}{section}

\def\openone

{\mathchoice

{\hbox{\upshape \small1\kern-3.3pt\normalsize1}}

{\hbox{\upshape \small1\kern-3.3pt\normalsize1}}

{\hbox{\upshape \tiny1\kern-2.3pt\SMALL1}}

{\hbox{\upshape \Tiny1\kern-2pt\tiny1}}}

\makeatletter

\newbox\ipbox

\newcommand{\ip}[2]{\left\langle #1\, , \,#2\right\rangle}
\newcommand{\diracb}[1]{\left\langle #1\mathrel{\mathchoice

{\setbox\ipbox=\hbox{$\displaystyle \left\langle\mathstrut
#1\right.$}

\vrule height\ht\ipbox width0.25pt depth\dp\ipbox}

{\setbox\ipbox=\hbox{$\textstyle \left\langle\mathstrut
#1\right.$}

\vrule height\ht\ipbox width0.25pt depth\dp\ipbox}

{\setbox\ipbox=\hbox{$\scriptstyle \left\langle\mathstrut
#1\right.$}

\vrule height\ht\ipbox width0.25pt depth\dp\ipbox}

{\setbox\ipbox=\hbox{$\scriptscriptstyle \left\langle\mathstrut
#1\right.$}

\vrule height\ht\ipbox width0.25pt depth\dp\ipbox}

}\right. }

\newcommand{\dirack}[1]{\left. \mathrel{\mathchoice

{\setbox\ipbox=\hbox{$\displaystyle \left.\mathstrut
#1\right\rangle$}

\vrule height\ht\ipbox width0.25pt depth\dp\ipbox}

{\setbox\ipbox=\hbox{$\textstyle \left.\mathstrut
#1\right\rangle$}

\vrule height\ht\ipbox width0.25pt depth\dp\ipbox}

{\setbox\ipbox=\hbox{$\scriptstyle \left.\mathstrut
#1\right\rangle$}

\vrule height\ht\ipbox width0.25pt depth\dp\ipbox}

{\setbox\ipbox=\hbox{$\scriptscriptstyle \left.\mathstrut
#1\right\rangle$}

\vrule height\ht\ipbox width0.25pt depth\dp\ipbox}

} #1\right\rangle}

\newcommand{\cj}[1]{\overline{#1}}

\newcommand{\bz}{\mathbb{Z}}

\newcommand{\br}{\mathbb{R}}

\newcommand{\bt}{\mathbb{T}}
\newcommand{\bn}{\mathbb{N}}

\def\blfootnote{\xdef\@thefnmark{}\@footnotetext}

\renewcommand{\mod}{\operatorname{mod}}

\input xy
\xyoption{all}
\usepackage{amssymb}

\newcommand{\supp}[1]{\text{supp} (#1)}

\def\-{^{-1}}

\begin{document}

\title[Quasiperiodic spectra]{QUASIPERIODIC SPECTRA AND ORTHOGONALITY FOR ITERATED FUNCTION SYSTEM
MEASURES}
\author{Dorin Ervin Dutkay}
\blfootnote{Research supported in part by a grant from the National Science Foundation DMS-0704191}
\address{[Dorin Ervin Dutkay] University of Central Florida\\
	Department of Mathematics\\
	4000 Central Florida Blvd.\\
	P.O. Box 161364\\
	Orlando, FL 32816-1364\\
U.S.A.\\} \email{ddutkay@mail.ucf.edu}

\author{Palle E.T. Jorgensen}
\address{[Palle E.T. Jorgensen]University of Iowa\\
Department of Mathematics\\
14 MacLean Hall\\
Iowa City, IA 52242-1419\\}\email{jorgen@math.uiowa.edu}
\thanks{} 
\subjclass[2000]{65F10, 42C15, 42A55, 42C05, 46E22, 28A80, 52C22, 52C23, 47A56, 65F30, 46F30}
\keywords{Iterated function systems, orthogonal bases, Hilbert space, Fourier series, Riesz products, special orthogonal functions, fractals, quasicrystals, aperiodic tilings, Hadamard matrix, matrix algorithms, nonlinear analysis.}

\begin{abstract}
  We extend classical basis constructions from Fourier analysis to attractors for affine iterated function systems (IFSs). This is of interest since these attractors have fractal features, e.g., measures with fractal scaling dimension. Moreover, the spectrum is then typically quasi-periodic, but non-periodic, i.e., the spectrum is a ``small perturbation'' of a lattice. Due to earlier research on IFSs, there are known results on certain classes of spectral duality-pairs, also called spectral pairs or spectral measures. It is known that some duality pairs are associated with complex Hadamard matrices. However, not all IFSs $X$ admit spectral duality. When $X$ is given, we identify geometric conditions on $X$ for the existence of a Fourier spectrum, serving as the second part in a spectral pair. We show how these spectral pairs compose, and we characterize the decompositions in terms of atoms. The decompositions refer to tensor product factorizations for associated complex Hadamard matrices.
\end{abstract}
\maketitle \tableofcontents

\section{Introduction}\label{intr}
  The idea of expanding $L^2$-functions on subsets $\Omega$ in Euclidean space into bases of more fundamental functions is central, and dates back to Fourier. It is of use in signal processing and in physics, but also of interest in its own right: Here we are thinking of Fourier series, orthogonal polynomials, eigenfunctions for Hamiltonians in physics, and wavelets; to mention only a few. These are instances of Marc Kac's question: ``Can you hear the shape of a drum?'' Each case suggests a natural choice of basis functions. We will consider a setting when the set $\Omega$ under consideration comes with some degree of selfsimilarity, and we will be asking for the possibility of choosing Fourier bases; i.e., we will examine the possibility of selecting orthonormal bases in $L^2(\Omega, \mu)$ where $\mu$ is a finite measure on $\Omega$ which reflects the intrinsic selfsimilarity. Such selfsimilarity arises for example in affine iterated function systems \cite{Hut81}, but it is much more general as we demonstrate. If $\Omega$ has non-empty interior, it is natural to take $\mu$ to be the restriction of Lebesgue measure. Hence we are faced with a pair of subsets in $\br^d$: (1) the set $\Omega$ itself, and (2) the points $\lambda$ which make up the frequencies in some candidate for a Fourier basis in $L^2(\Omega)$.

    In the discussion below, we recall the history of the problem, and earlier results by a number of authors which are relevant for our present work. The central question we address here is this: ``To what extent may some given set $\Omega$ in d dimensions be built up from atoms of fundamental blocks in such a way that the spectral data for the ``atoms'' determine that of $\Omega$ itself?'' Even if the spectral data for the atoms is periodic, we show that for composite systems, the expectation is quasiperiodicity in a sense we make precise in section 3 below.

    Our work is inspired by \cite{Fug74, IKT01, Lo67, Lab02} among others.

We consider open subsets $\Omega$ in $\br^d$ of finite positive Lebesgue measure. Our focus is on the case when the Hilbert space $L^2(\Omega)$ has an orthogonal Fourier basis, i.e., an orthogonal basis complex exponentials. The measure on $\Omega$ is taken to be the restriction of $d$-dimensional Lebesgue measure. The exponents in such an orthogonal basis will then form a discrete subset $\Lambda$ in $\br^d$. We say that $(\Omega, \Lambda)$ is a spectral pair and $\Omega$ is a spectral set.

     We identify a geometric condition which characterizes spectral pairs arising as attractors of iterated function systems (IFSs), i.e., from a finite set of affine mappings in $\br^d$.

 We analyze sets of the form $A+[0,1]$ where $A$ is some finite set of integers, and find conditions when such a set is spectral (Theorem \ref{thspec}). We characterize those sets which are attractors of an affine IFS (Theorem \ref{thtile} and show that they are spectral sets (Theorem \ref{thfinsp}). We construct a new class of spectral measures (Theorem \ref{thpro}), and obtain a counterexample to a conjecture of \L aba and Wang (Example \ref{ex3_6}). We present an example of a measure which has an infinite family of mutually orthogonal exponentials but is not spectral (Proposition \ref{pr3_24}). We show how new spectra can be constructed from old for some fractal measures (Lemma \ref{lemnews} and Theorem \ref{thinf}). We construct a {\it connected} spectral domain in $\br^3$ which does not tile $\br^3$ by any lattice (Example \ref{exst}).

   We introduce more general spectral pairs than the $(\Omega,\Lambda)$ systems, including a pairing for finite subsets in $\br^d$, and from IFSs. And we introduce an operation on spectral pairs. Our idea is to
identify an interplay between finite spectral pairs on the one hand, and a
class of infinite Euclidean ones on the other, those built on affine
iterated function system (IFS) measures, see Definition \ref{defaifs}. With tools from
IFS-theory, this then allows us to exploit our new results on finite systems
in extending some of the classical constructions from Fuglede's paper
\cite{Fug74}.

Section 3 contains several new results: (a) A FFT-type algorithm (Corollary \ref{cor3_16}) in 1D of building molecules of spectral pairs $(\Omega, \Lambda)$ from atoms. (b) For this class of spectral pairs $(\Omega, \Lambda)$, when $\Omega$ is fixed, we find all the possible sets $\Lambda$ which serve as spectra (Theorem \ref{thspec}.) In section 4 we consider systems in higher dimensions, with special attention to the case when $\Omega$ is both open and connected.
    
             The broader motivation for our paper is a set of intriguing connections between tiles, spectrum and wavelet analysis. To a large degree, the role of scaling operators has been missing in many early approaches to spectral-tile duality. The advent of wavelets \cite{Dau92} did much to remedy this. Some early papers stressing the role played by scaling and selfsimilarity in spectrum-tile duality and in wavelets are  \cite{Law91, BrJo99, JoPe99, BJR99}, and especially \cite{GrMa92} which make useful connections to signal processing in engineering. Our main results concern spectral properties implied by selfsimilarity.

The implications of this selfsimilarity (i.e., similarity up to a
suitable scaling operation, or a group of affine mappings) take several
forms:
   Our Corollary \ref{corr} below identifies the transformation rules for the
action of the affine group $\mathcal A_d$ in $\br^d$ on the finite Borel measures on $\br^d$,
and on the subsets $\Lambda$ in $\br^d$ which can occur as spectra of these
measures.

   An affine IFS (\cite{Hut81}, Theorem \ref{thhut}) and its invariant measures $\mu$ are
defined from a prescribed finite subset $\mathcal F$ in the group $\mathcal A_d$. Different
choices of subsets $\mathcal F$ in $\mathcal A_d$ yield different affine IFSs. In our separate
results Theorems \ref{thpro}, \ref{thspec}, \ref{thifs}, \ref{thfinsp} and \ref{thinf} we derive detailed
spectral data for affine IFSs. Specifically, we derive quasiperiodic
spectral properties of these IFS-invariant measures $\mu$ (see equation \eqref{eqmub} below)
making use the scaling-similarity implied by $\mathcal F$-invariance; i.e., spectral
data for measures defined from an $\mathcal F$-invariance property for a prescribed
finite subset $\mathcal F$ in the group $\mathcal A_d$.

 The introduction of a suitable scaling operation further makes a connection to
selfsimilar structures that arise in electrical networks (e.g., \cite{Pow76}) and for the
IFS-fractals $(X,\mu)$ of Kigami, Lapidus and Strichartz; i.e., Sierpinski gaskets, Sierpinski
carpets etc, see e. g., \cite{KiLa01}, \cite{St06a}, and \cite{St06b}. With the notation $(X,\mu)$, it is understood
that for a particular IFS, the associated measure $\mu$ is a Hutchinson equilibrium measure \cite{Hut81}
with support $X$. The set $X$ may be a Cantor set, or a  Sierpinski gasket. Both the measure and its support
have a scaling dimension $\delta$, typically a fraction. For the Cantor set it is $\delta =\log_3(2)$.

      What sets the fractal IFSs apart from the original Euclidean systems $(\Omega, \Lambda)$
without scale-similarity (see \cite{Fug74}) is that in the fractal case, the orthogonal Fourier functions 
$\Lambda$ in $L^2(X,\mu)$ form {\it local bases} when the measure $\mu$ in question is a fractal IFS equilibrium measure
(Definition \ref{defaifs}). More precisely, Strichartz \cite{St06b} proved that this local feature (shared by wavelet bases) accounts for better approximations, i.e., for better convergence of the associated $\Lambda$-Fourier series. Specifically \cite{St06b}, in the IFS case the $\Lambda$-Fourier series converges for all continuous functions on $X$.

                A source of inspiration for our results is an idea in a recent sequence of papers by \L aba, and by \L aba and Wang \cite{LaWa06, Lab01, Lab02, LaWa02}, as well as \cite{JoPe99}, results which suggest the usefulness in studying spectral tile duality with the aid of ``fundamental building blocks'';  see the next two sections below for details.

      While the original Fuglede conjecture \cite{Fug74} is known to be negative for Lebesgue measure restricted to subsets in $\br^d$ when $d$ is 3 or more (this is work beginning with \cite{Tao04} then \cite{FMM06,KM06}), so far little is known in the way of complete spectral/tile results for small $d$, even for $d = 1$. And if the measures $\mu$ under consideration arise as equilibrium measures for affine iterated function systems (IFS), as is often the case for fractal measures, again then there are only partial results in the literature regarding connections between geometry and spectra (in the form of orthogonal Fourier bases in $L^2(\mu)$.)

      One of the conclusions from our present work is that a rich class of IFS-measures may have the form $\mu =$ Lebesgue measure restricted to a suitably chosen subset in $\br^d$ of finite positive (Lebesgue) measure, i.e., measures arising from restriction to finite geometries in $\br^d$. Moreover we show that this geometry for configurations in $\br^d$ is closely connected to the question of when $\mu$ is a spectral measure. Since spectral results in low dimensions are sparse, we feel that a closer examination of such new approaches is worthwhile.

\section{Definitions}\label{defi}

 In this section, we identify an interplay between finite spectral pairs
with those built on infinite iterated function system (IFS) measures. This
allows us in section 3 to resolve a conjecture of \L aba-Wang, see Conjecture \ref{conlw}
.

\begin{definition}\label{defspec}
For $\lambda\in\br^d$, let $e_\lambda(x):=e^{2\pi i\lambda\cdot x}$, $x\in\br^d$.

We say that a measure $\mu$ on $\br^d$ is a {\it spectral measure} if there exists a set $\Lambda\subset\br^d$, such that $\{e_\lambda\,|\,\lambda\in\Lambda\}$ is an orthogonal basis in $L^2(\mu)$. In this case $\Lambda$ is called a {\it spectrum} for the measure $\mu$. We call $(\mu,\Lambda)$ a {\it spectral pair}.

We say that a finite set $A\subset\br^d$ is {\it spectral} if the atomic measure $\delta_A:=\frac{1}{\#A}\sum_{a\in A}\delta_a$ is spectral. $\#A$ denotes the cardinality of $A$, and $\delta_a$ is the Dirac measure at $a$. A set $\Lambda$ is called a {\it spectrum} for $A$ if it is a spectrum for $\delta_A$. We call $(A,\Lambda)$ a {\it spectral pair}.

We say that a Lebesgue measurable set $\Omega$ of positive finite Lebesgue measure in $\br^d$ is {\it spectral} if the Lebesgue measure restricted to $\Omega$ is spectral. A {\it spectrum} $\Lambda$ for $\Omega$ is any spectrum for the Lebesgue measure on $\Omega$. $(\Omega,\Lambda)$ is called a {\it spectral pair}.
\end{definition}

\begin{definition}\label{defsp}
Let $G$ be an abelian group. We say that a subset $A$ of $G$ {\it tiles $G$} if there exists a set $\mathcal T$ such that $(A+t)_{t\in\mathcal T}$ is a partition of $G$ up to Haar measure zero, i.e., if $\mu_G$ is the Haar measure on $G$, then
$$\mu_G\left(G\setminus\bigcup_{t\in\mathcal T}(A+t)\right)=0,\quad
\mu_G\left((A+t)\cap(A+t')\right)=0,\quad\mbox{ for all }t,t'\in \mathcal T, t\neq t'.$$

We call $\mathcal T$ a {\it tile set} for $A$. We say that {\it $A$ tiles $G$ with $\mathcal T$}.

When $G$ is a discrete group and $A$ tiles $G$ with $\mathcal T$, we use the notation $A\oplus\mathcal T=G$.

If $A$,$B$ are subsets of a discrete group, we use the notation $A\oplus B=C$, if every $c$ in $C$ can be written uniquely as $c=a+b$ with $a\in A$ and $b\in B$.
\end{definition}

\begin{definition}
If $\mu$ is a finite measure on $\br^d$, then we denote by $\hat\mu$ the {\it Fourier transform} of $\mu$:
$$\hat\mu(x)=\int_{\br^d}e^{2\pi it\cdot x}\,d\mu(t)\quad(x\in\br^d).$$
\end{definition}

Note that if $\ip{\cdot}{\cdot}$ denotes the inner product in $L^2(\mu)$, then 
$$\ip{e_\lambda}{e_{\lambda'}}=\hat\mu(\lambda-\lambda').$$
   As motivation for our problem, we recall the statement of the Fuglede conjecture, although it is now known to be negative in general: The question, or conjecture, from \cite{Fug74} was, if for a given measurable subset $\Omega$ in $\br^d$ of finite positive Lebesgue measure, the following two conditions \eqref{equnu} and \eqref{eqdoi} are equivalent:

\begin{equation}\label{equnu}
\Omega\mbox{ tiles }\br^d\mbox{ with translations.}
\end{equation}
\begin{equation}\label{eqdoi}
\Omega\mbox{ has a spectrum, i.e., }L^2(\Omega)\mbox{ has an ONB of Fourier frequencies.}
\end{equation}

     This conjecture is known to be negative in $\br^d$ for $d \geq 4$; see e.g., \cite{Tao04}, and now also disproved in 3 dimensions, see \cite{FMM06} and \cite{KM06}, but it is open in lower dimensions.

     Hence in the literature, starting with \cite{PW01, LW96}, a number of authors have placed additional conditions on the sets in \eqref{equnu} and the spectra in \eqref{eqdoi} with view to more definite results. For example, in \cite{Ped04a, Ped04b} Pedersen introduced an intriguing ``dual spectral-set-conjecture''.

       Here we address the question for $d = 1$, of whether a tiling property for $\Omega$ together with a degree of selfsimilarity (details below) implies the spectral property. Even though this is then more restrictive, more specific, it is of interest even for dimension $d = 1$.

\begin{definition}\label{defaifs}
Let $A$ be a $d\times d$ expansive integer matrix. We say that a matrix is {\it expansive} if all its eigenvalues have absolute value $>1$. Let $B$ be a finite subset of $\br^d$. We call the family of maps $(\tau_b)_{b\in B}$, 
$$\tau_b(x)=A^{-1}(x+b),\quad(x\in\br^d,b\in B),$$
an {\it affine iterated function system (affine IFS)}. 
\end{definition}

\begin{theorem}\label{thhut}\cite{Hut81}
Let $(\tau_b)_{b\in B}$ be an affine IFS. There is a unique compact subset $X_B$ of $\br^d$ such that
\begin{equation}\label{eqxb}
X_B=\bigcup_{b\in B}\tau_b(X_B)
\end{equation}
There exists a unique probability measure $\mu_B$ on $\br^d$ that satisfies the following {\it invariance equation}
\begin{equation}\label{eqmub}
\int_{\br^d}f\,d\mu_B=\frac{1}{\#B}\sum_{b\in B}\int_{b\in B}f\circ\tau_b\,d\mu_B,\quad(f\in C_c(\br^d)).
\end{equation}
Moreover $\mu_B$ is supported on $X_B$.
\end{theorem}

\begin{definition}
The compact set $X_B$ is called the {\it attractor} of the affine IFS $(\tau_b)_{b\in B}$. The measure $\mu_B$ is called the {\it invariant measure} of the IFS $(\tau_b)_{b\in B}$.
\end{definition}

\begin{theorem}\label{thdj}\cite{DuJo06} One dimension.
Let $\tau_b(x)=A^{-1}(x+b)$, $b\in B$, $x\in\br$, be an affine IFS, with $A\in\bz$, $A\geq 2$ and $B\subset\bz$, $0\in B$. Assume that $B$ is spectral with spectrum $\frac{1}{A}L$ for some subset $L$ of $\bz$, with $0\in L$. Then the invariant measure $\mu_B$ is a spectral measure. 
If in addition $\gcd(B)=1$, then the measure $\mu_B$ has a spectrum contained in $\bz$.
\end{theorem}

\begin{proof}
The fact that $\mu_B$ is spectral is proved in \cite{DuJo06}.
We only need to prove the last statement, that the spectrum constructed in \cite{DuJo06} is contained in $\bz$. Recall that this spectrum is the smallest set $\Lambda$ which contains $-C$ for all $\hat\delta_B$-cycles $C$, and such that $A\Lambda+L\subset\Lambda$. For the definition of $\hat\delta_B$-cycles we refer to \cite{DuJo06}. A point $c$ in a $\hat\delta_B$-cycle has the property $|\hat\delta_B(c)|=1$. But this implies that 
$$\frac1{\#B}\left|\sum_{b\in B}e^{2\pi i bc}\right|=1.$$
Using the triangle inequality, and since $0\in B$, we see that all the terms in the sum must be equal to $1$. So 
$bc\in\bz$ for all $b\in B$. Since $\gcd(B)=1$, there exist integers $m_b$, for all $b\in B$, such that $\sum_{b\in B}m_bb=1$
then $c=\sum_{b\in b}m_bbc\in\bz$. Therefore $-C$ is contained in $\bz$ for all $\hat\delta_B$-cycles, so the smallest set $\Lambda$ that contains $-C$ and satisfies $A\Lambda+L\subset\Lambda$ is contained in $\bz$. 
\end{proof}

\section{Spectral theory for measures}\label{spectral}

  In the theory of quasi crystals in higher dimensions $d$ (see especially \cite{BaMo00, BaMo01}), one often encounters finite measures $\mu$ on $\br^d$ with spectrum consisting of discrete subsets $\Lambda$ which posses a certain quasi-periodicity. While the interesting sets $\Lambda$ are not rank-$d$ lattices, they are in a certain sense ``small perturbations'' of lattices.

       Our first lemma is a characterization of a pair $(\mu, \Lambda)$ in $\br^d$ when $\mu$ is a finite measure and $\Lambda$ is a subset of $\br^d$. It is a necessary and sufficient condition for $(\mu, \Lambda)$ to be a spectral pair. While it was noticed also in \cite{JoPe99}, we sketch the details below for the benefit of the reader.

       \begin{lemma}\label{lemspectral}
Let $\mu$ be a probability measure on $\br^d$. Then $\mu$ has spectrum $\Lambda$ if and only if
\begin{equation}\label{eqspectrum}
1=\sum_{\lambda\in\Lambda}|\hat\mu(x+\lambda)|^2,\quad(x\in\br^d).
\end{equation}
\end{lemma}

\begin{proof}
We have for all $x,y\in\br^d$:
$$\ip{e_x}{e_y}=\int e^{2\pi i(x-y)\cdot t}\,d\mu(t)=\hat\mu(x-y).$$
Therefore the necessity of \eqref{eqspectrum} follows from the Parseval equality.

Conversely, if \eqref{eqspectrum} is satisfied, then take $x=-\lambda'$ for some $\lambda'\in\Lambda$. Since $\mu$ is a probability measure, $\hat\mu(0)=1$, and equation \eqref{eqspectrum} implies that $\hat\mu(\lambda-\lambda')=0$ for all $\lambda\in\Lambda$, $\lambda\neq\lambda'$. Thus $(e_\lambda)_{\lambda\in\Lambda}$ forms an orthonormal family of vectors in $L^2(\mu)$. We have to check only that it is complete.

Let $H$ be the closed span of the functions $(e_\lambda)_{\lambda\in\Lambda}$ and let $P$ the projection onto $H$. Since $(e_\lambda)_{\lambda\in\Lambda}$ is an orthonormal basis for $H$, we have for all $x\in\br^d$:
$$\|Pe_{-x}\|^2=\sum_{\lambda\in\Lambda}|\ip{e_{\lambda}}{e_{-x}}|^2=\sum_{\lambda\in\Lambda}|\hat\mu(x+\lambda)|^2=1=\|e_{-x}\|^2.$$
But this implies that $e_{-x}$ is in $H$. Since $x$ is arbitrary, we can use the Stone-Weierstrass theorem to conclude that $H=L^2(\br^d)$.
\end{proof}

\begin{lemma}\label{lemftinv}
Let $(\tau_b)_{b\in B}$ be an affine IFS, and let $\mu_B$ be its invariant measure. Then
$$\hat\mu_B(A^Tx)=\hat\delta_B(x)\hat\mu_B(x),\quad(x\in\br^d).$$

$$\hat\mu_B(x)=\prod_{n=1}^\infty\hat\delta_B((A^T)^{-n}x),\quad(x\in\br^d).$$
\end{lemma}
\begin{proof}
Just take the Fourier transform of the invariance equation in Theorem \ref{thhut}. Iterating the scaling equation and using the fact that $A$ is expansive, one gets the infinite product formula. The product is uniformly convergent on compact subsets. 
\end{proof}

\subsection{Action of the affine group}

 In this section we prove a rigidity theorem for the action of the affine group $\mathcal A_d$ in $\br^d$ on the finite Borel measures on $\br^d$, and on the subsets $\Lambda$ in $\br^d$ which can occur as spectra of these measures.

\begin{definition}
Let $d\in\bn$ be given. Set 
$$\mathcal M_d:=\{\mbox{ all finite positive Borel measures on }\br^d\},\quad \mathcal{SM}_d:=\{\mu\in\mathcal M_d\,|\,\mu\mbox{ is spectral }\}.$$
$\mathcal A_d$:=the affine group of $\br^d$, i.e., all invertible affine transformations $x\mapsto Vx+s: \br^d\rightarrow\br^d$, where $V\in GL_d$, and $s\in\br^d$.

For $\mu\in\mathcal M_d$, $a\in\mathcal A_d$, set
\begin{equation}\label{eqr1}
R(a)\mu:=\mu\circ a^{-1}
\end{equation}

If $\mu\in\mathcal SM_d$, and if $(\mu,\Lambda)$ is a spectral pair, we say that $\Lambda\in\mathcal S(\mu)$.
\end{definition}

\begin{corollary}\label{corr}
Via the formula \eqref{eqr1}, the affine group $\mathcal A_d$ acts as a transformation group on $\mathcal M_d$, and $\mathcal{SM}_d$ is invariant, i.e., 
\begin{equation}\label{eqr2}
\mbox{ If }a\in\mathcal A_d,\mu\in\mathcal{SM}_d\mbox{ then }R(a)\mu\in\mathcal{SM}_d.
\end{equation}

If $a\in\mathcal A_d$, and $a(x)=Vx+s$, then 
\begin{equation}\label{eqr3}
\Lambda\in\mathcal S(R(a)\mu)\mbox{ iff }V^T\Lambda\in\mathcal S(\mu).
\end{equation}

\end{corollary}
\begin{example}
Let $\mu$ be the IFS measure in $\br$ given by $x\mapsto\frac x4$, and $x\mapsto\frac{x+2}4$, see Definition \ref{defaifs} and Theorem \ref{thhut}, i.e., $A=4$, and $B=\{0,2\}$. It is known \cite{JoPe98} that $\mu\in\mathcal{SM}_1$. It follows that the IFS measure $\mu_1$ for $A=4$ and $B_1=\{0,1\}$ is in $\mathcal{SM}_1$, and that 
$$\Lambda\in\mathcal S(\mu)\mbox{ iff } 2\Lambda\in\mathcal S(\mu_1).$$
\end{example}

\begin{proof}{\it of Corollary \ref{corr}.}

Let $a\in\mathcal A_d$. Then there are $V\in GL_d$, and $s\in\br^d$, such that $a(x)=Vx+s$. Consider $\mu\in\mathcal M_d$, and the Fourier transform $\widehat{R(a)\mu}$. We then have
\begin{equation}\label{eqr4}
\widehat{R(a)\mu}(t)=\int e_t(x)\,d(\mu\circ a^{-1})(x)=\int e_t(Vx+s)\,d\mu(x)=e_t(s)\hat\mu(V^Tt).
\end{equation}

If $\mu\in\mathcal{SM}_d$ and $\Lambda\in\mathcal S(\mu)$, we then get the following identity for $\widehat{R(a)\mu}$, where the notation $W:=V^T$ is used:
$$\sum_{\lambda\in\Lambda}|\widehat{R(a)\mu}(t+W^{-1}\lambda)|^2=\sum_{\lambda\in\Lambda}|\hat\mu(W(t+W^{-1}\lambda))|^2=\sum_{\lambda\in\Lambda}|\hat\mu(Wt+\lambda)|^2=1,$$
for all $t\in\br^d$ by Lemma \ref{lemspectral}. As a result we conclude that $W^{-1}\Lambda\in\mathcal S(R(a)\mu)$ as claimed.

\end{proof}

\subsection{Induction from finite measures}
The next lemma offers a complete characterization of finite spectral
pairs, and it will be needed in the proof of our results on building new
spectral pairs from ``old ones''.

 In our separate results Theorems \ref{thpro} and \ref{thspec} below we combine induction from finite measures with an analysis of Hadamard matrices in deriving detailed spectral data for affine IFSs. Specifically, we establish quasiperiodic spectral properties of invariant IFS-measures (see equation \eqref{eqmub}) making use of the scaling-similarity.

\begin{lemma}\label{lemfini}
Let $F:=\{x_1,\dots,x_p\}$ be some a finite set of distinct points in $\br^d$. Let $\delta_{F}$ be the measure
$$\delta_F=\frac1p\sum_{i=1}^p\delta_{x_i}.$$
\begin{enumerate}
\item
The set $F$ has spectrum $\Lambda$ if and only if $\#\Lambda=p$, $\Lambda=\{\lambda_1,\dots,\lambda_p\}$, and the matrix
$$\frac{1}{\sqrt p}\left(e^{2\pi i x_i\cdot\lambda_j}\right)_{i,j=1}^p$$
is unitary.
\item The Fourier transform of $\delta_F$ is 
$$\hat\delta_F(t)=\frac{1}{p}\sum_{i=1}^pe^{2\pi ix_i\cdot t},\quad(t\in\br^d).$$
\item The set $F$ has spectrum $\Lambda$ if and only if 
$$\sum_{\lambda\in\Lambda}|\hat\delta_F(x+\lambda)|^2=1,\quad(x\in\br^d).$$
\end{enumerate}
\end{lemma}
\begin{proof}
The Hilbert space $L^2(\delta_F)$ clearly has dimension $p$. Therefore any spectrum for $\delta_F$ will have cardinality $p$. Also
$$\ip{e_\lambda}{e_\lambda'}_{L^2(\delta_F)}=\frac{1}{p}\sum_{i=1}^p e^{2\pi ix_i\cdot(\lambda-\lambda')}=\delta_{\lambda,\lambda'}.$$
This translates into the rows of the matrix being orthonormal.

(ii) follows by direct computation.

(iii) follows from Lemma \ref{lemspectral}.
\end{proof}

\begin{lemma}\label{lem2}
Let $\mu_1,\mu_2$ be two probability measures on $\br^d$. Suppose the following assumptions hold:
\begin{enumerate}
\item $\mu_1$ and $\mu_2$ are spectral measures with spectra $\Lambda_1$ and $\Lambda_2$ respectively.
\item $\hat\mu_2(x+\lambda_1)=\hat\mu_2(x)$ for all $x\in\br^d$ and $\lambda_1\in\Lambda_1$.
\end{enumerate}

Then the convolution measure $\mu_1\ast\mu_2$ is a spectral measure with spectrum $\Lambda_1\oplus\Lambda_2$.

\end{lemma}
\begin{proof}
We use Lemma \ref{lemspectral}. 

We have
$$\sum_{\lambda_1\in\Lambda_1,\lambda_2\in\Lambda_2}|\widehat{\mu_1\ast\mu_2}(x+\lambda_1+\lambda_2)|^2=\sum_{\lambda_2\in\Lambda_2}\sum_{\lambda_1\in\Lambda_1}|\hat\mu_1(x+\lambda_1+\lambda_2)|^2|\hat\mu_2(x+\lambda_1+\lambda_2)|^2=$$
$$=\sum_{\lambda_2\in\Lambda_2}|\hat\mu_2(x+\lambda_2)|^2\sum_{\lambda_1\in\Lambda_1}|\hat\mu_1(x+\lambda_1+\lambda_2)|^2=\sum_{\lambda_2\in\Lambda_2}|\hat\mu_2(x+\lambda_2)|^2=1.$$
First, this proves that an element $\lambda\in\Lambda_1+\Lambda_2$ can be written uniquely as $\lambda=\lambda_1+\lambda_2$ with $\lambda_1\in\Lambda_1$ and $\lambda_2\in\Lambda_2$. Otherwise, take $x=-\lambda$, and the lefthand side of the equality is greater than $2$, because $\widehat{\mu_1\ast\mu_2}(0)=1$.

This, and Lemma \ref{lemspectral} proves that $\mu_1*\mu_2$ has spectrum $\Lambda_1\oplus\Lambda_2$.
\end{proof}

\begin{corollary}\label{corfinite}
Consider the following assumptions
\begin{enumerate}
\item
The probability measure $\mu_1$ on $\br$ has spectrum $\Lambda_1$ contained in $\bz$. 
\item $F$ is a finite subset of $\bz$ with spectrum $\Lambda_2$.
\end{enumerate}
Then $\mu_1\ast\delta_F$ has spectrum $\Lambda_1\oplus\Lambda_2$.
\end{corollary}

\begin{example}\label{ex3_6}
Let $\mu_1$ be the Lebesgue measure on $[0,1]$. Let $\mu_2=\delta_{\{0,2\}}$. The convolution $\mu_1\ast\mu_2$ is the Lebesgue measure on $[0,1]+\{0,2\}=[0,1]\cup[2,3]$, renormalized so that it is a probability measure. This can be seen from the following
\begin{lemma}
Let $\mu$ be a measure on $\br^d$, and $F$ a finite subset of $\br^d$. Then for all $f\in C_c(\br^d)$,
$$\int f\,d\mu\ast\delta_F=\frac1{\#F}\sum_{\alpha\in F}\int f(x+\alpha)\,d\mu(x).$$
\end{lemma}
The proof requires a one line computation.

We have that $\mu_1$ is spectral with spectrum $\Lambda_1=\bz$. Also, with Lemma \ref{lemfini}, $\mu_2$ is spectral with spectrum $\Lambda_2:=\{0,\frac14\}$. We have 
$$\hat\mu_2(t)=\frac{1}{2}(1+e^{2\pi i2t}),\quad(t\in\bt).$$
It is clear that all the assumptions of Lemma \ref{lem2} are satisfied. Therefore we obtain that the renormalized Lebesgue measure on $[0,1]\cup[2,3]$ has spectrum $\bz+\{0,\frac14\}$.  No lattice $\Lambda$ is a spectrum for $[0,1]\cup[2,3]$ because this set does not tile $\br$ by any lattice.
\end{example}

\begin{remark}\label{remlw}
In \cite{LaWa02}, \L aba and Wang proposed the following conjecture.

\begin{conjecture}\cite{LaWa02}\label{conlw}
Let $\mu$ be the invariant measure associated with the IFS $\phi_j(x)=\rho(x+a_j)$, $1\leq j\leq q$, with probability weights $p_1,\dots,p_q>0$, where $|\rho|<1$. Suppose that $\mu$ is a spectral measure. Then
\begin{enumerate}
\item[(a)] $\rho=\frac1N$ for some $N\in\bz$.
\item[(b)] $p_1=\dots=p_q=\frac1q$.
\item[(c)] Suppose that $0\in\mathcal A=\{a_j\}$. Then $\mathcal A=\alpha\mathcal D$ for some $\alpha\in\br$ and $\mathcal D\subset\bz$. Furthermore, $\mathcal D$ must be a complementing set $(\mod N)$, i.e., there exists a set $\mathcal E\subset \bz$ such that $\mathcal D\oplus \mathcal E$ is a complete residue system $(\mod N)$. 
\end{enumerate}
\end{conjecture}

The set in Example \ref{ex3_6} provides a counterexample to the last statement in this conjecture. Indeed, the set 
$[0,1]\cup[2,3]$ is the attractor of the IFS $\tau_b(x)=\frac{1}{4}(x+b)$ with $b\in B:=\{0,1,8,9\}$. The invariant measure of this IFS is the normalized Lebesgue measure on $[0,1]\cup[2,3]$. This can be proved by checking the invariance equations in Theorem \ref{thhut}. Since $0\equiv 8\mod 4$, the set $B$ is not complementing.

In Theorem \ref{thpro} we will provide a larger class of affine iterated function systems which yield examples of spectral measures which contradict the \L aba-Wang conjecture.

\end{remark}
\begin{example}\label{exjp}
Let $\mu_4$ be the invariant measure of the affine IFS $\tau_b(x)=A^{-1}(x+b)$ with $A=4$ and $B=\{0,2\}$. It was proved in \cite{JoPe98} that this is a spectral measure with spectrum $\{\sum_{k=0}^n4^kl_k\,|\,k\in\{0,1\}\}$.

The measure $\mu_4\ast\delta_{\{0,2\}}$ has spectrum
$$\left\{\sum_{k=0}^n4^kl_k\,|\,l_k\in\{0,1\}\right\}+\{0,\frac14\}.$$
\end{example}

\begin{example}
$[0,2]\cup[5,6]$ has spectrum $\bz+\{0,\frac13,\frac23\}=\frac13\bz$, a lattice. To see this take $F=\{0,1,5\}$ and $\Lambda_2=\{0,\frac13,\frac23\}$.
\end{example}

\begin{example}
Let $d$ be an integer $d\neq0$. Then $[0,1]\cup[d,d+1]$ has spectrum $\bz+\{0,\frac1{2d}\}$. Take $F=\{0,d\}$ and $\Lambda_2=\{0,\frac1{2d}\}$.
\end{example}

\begin{corollary}\label{cor3_11}
Let $D$ be a finite subset of $\br$. Assume the following conditions are satisfied:
\begin{enumerate}
\item[(a)] $D=\oplus_{k=1}^nb_kC_k$ with $b_k\in\br\setminus\{0\}$ and $C_k\subset\bz$.
\item[(b)] For all $k\in\{1,\dots,n\}$, the set $C_k$ has spectrum $\frac{1}{a_k}L_k$ with $a_k\in\br\setminus\{0\}$ and $L_k\subset\bz$.
\item[(c)] $\frac{b_{k+1}}{a_jb_j}\in\bz$ for all $k,j\in\{1,\dots,n-1\}$, $j\leq k$.
\end{enumerate}
Then $D$ has spectrum $\oplus_{k=1}^n\frac{1}{a_kb_k}L_k$.
\end{corollary}

\begin{proof}
We will use Lemma \ref{lem2}, and the fact that for two finite subsets $A,B$ of $\br$ qith $A\oplus B=C$, one has $\delta_{A\oplus B}=\delta_A\ast\delta_B$.
We prove the corollary by induction. For $j=1$, $b_1C_1$ has spectrum $\frac{1}{a_1b_1}L_1$.

Assume by induction that for a $j<n$, the set $D_j:=\oplus_{k=1}^jb_kC_k$ has spectrum $S_j:=\oplus_{k=1}^j\frac{1}{a_kb_k}L_k$. We have to check condition (ii) in Lemma \ref{lem2}.

For any $l_k\in L_k$, $k=1,j$, 
$$\hat\delta_{b_{j+1}C_{j+1}}(x+\sum_{k=1}^j\frac{1}{a_kb_k}l_k)=\frac{1}{\#L_j}\sum_{c_{j+1}\in C_{j+1}}e^{2\pi ib_{j+1}c_{j+1}(x+\sum_{k=1}^j\frac{1}{a_kb_k}l_k)}=\hat\delta_{b_{j+1}C_{j+1}}(x),$$
and, in the last equality, we used condition (c) and $C_{j+1}, L_k\subset\bz$.

Then, with Lemma \ref{lem2}, $\oplus_{k=1}^{j+1}b_kC_k$ has spectrum $\oplus_{k=1}^{j+1}\frac{1}{a_kb_k}L_k$. 
The corollary follows by induction.
\end{proof}

\begin{corollary}\label{cor3_12}
Let $D$ be a finite subset of $\br$. Suppose the following conditions are satisfied.
\begin{enumerate}
\item[(a)] $D=\oplus_{k=1}^n a^{p_k}C_k$, where $0\leq p_1<p_2<\dots<p_n$ are some integers, $a\in\bz, a\geq1$, and $C_k\subset\bz$ for all $k\in\{1,\dots,n\}$.
\item[(b)] For all $k\in\{1,\dots,n\}$ the set $C_k$ has spectrum $\frac{1}{a}L_k$ for some $L_k\subset\bz$.
\end{enumerate}
Then the set $D$ has spectrum $\oplus_{k=1}^n\frac{1}{a^{p_k+1}}L_k$.

\end{corollary}

 The next lemma shows that there is an implication which is converse to
that of Lemma \ref{lem2}.

    It applies to a general class of affine iterated function systems
in $\br^d$. The next lemma offers a condition for when IFS-measures may
be factored as convolutions of more basic building block, where this
convolution-factorization is understood in the sense of Lemma \ref{lem2}.

    The general setup is as follows: We consider an initial IFS in $\br^d$
defined from a given and fixed $d\times d$ matrix $A$ and a finite subset $B$ in $\br^d$.
The pair $(A, B)$ determines an invariant measure $\mu = \mu_{A,B}$.

   We are assuming that the matrix $A$ is a $p$-fold product, i.e., $A = a^p$ for
some other $d\times d$ matrix $a$; and moreover that there is a compatible additive
decomposition of the set $B$. Geometrically, the pair $(A,B)$ factors into a
composition of ``atoms''.

   Under these conditions on the pair $(A,B)$, we obtain a convolution
factorization for the measure $\mu = \mu_{A,B}$ in terms of the atoms.

\begin{lemma}\label{lemdecomp}
Let $(\tau_b)_{b\in B}$ be an affine IFS defined by a $d\times d$
integer matrix $A$ and a finite subset $B$ of $\bz^d$. Consider the
following assumptions:
\begin{enumerate}
\item[(a)]  $A=a^p$ for some matrix $a$ and some $p\geq2$.
\item[(b)] $B=a^{n_0p}C_0\oplus a^{n_1p+1}C_1\oplus\dots\oplus
a^{n_{p-1}p+p-1}C_{p-1}$ for some $n_0,\dots,n_{p-1}\geq 0$ and some
finite subsets $C_0,\dots C_{p-1}$ of $\bz^d$ with the property that
if $c,c'\in C_i$ and $c-c'\in a\bz^d$ then $c=c'$.
\end{enumerate}
Then
$$\mu_{A,B}=\mu_{a^p,C_0\oplus aC_1\oplus\dots\oplus a^{p-1}C_{p-1}} *
\delta_F,$$ where
\begin{equation}\label{eqf}
F=\oplus_{k=0}^{p-1}\oplus_{l=0}^{n_k-1}a^{lp+k}C_k.
\end{equation}

If in addition $C_0=C_1=\dots =C_{p-1}=:C$, then $\mu_{a^p,C_0\oplus
aC_1\oplus\dots\oplus a^{p-1} C_{p-1}}=\mu_{a, C}$, so
$\mu_{A,B}=\mu_{a,C} *\delta_F$.
\end{lemma}

\begin{proof}
$$\hat\delta_B(x)=\frac{1}{\#B}\sum_{b\in B}e^{2\pi i b\cdot x}=
\frac{1}{\#B}\sum_{c_0\in C_0,\dots c_{p-1}\in C_{p-1}}e^{2\pi
i(a^{n_0p}c_0+a^{n_1p+1}c_1+\dots+a^{n_{p-1}p+p-1}c_{p-1})\cdot x}=$$
$$\prod_{k=0}^{p-1}\frac{1}{\#C_k}\sum_{c_k\in C_k}e^{2\pi i c_k\cdot
(a^T)^{n_kp+k}x}=\prod_{k=0}^{p-1}\hat\delta_{C_k}((a^T)^{n_kp+k}x).$$
Then

$$\hat\mu_{A,B}(x)=\prod_{n=1}^\infty\hat\delta_B((a^T)^{-np}x)=\prod_{n=1}^\infty\prod_{k=0}^{p-1}\hat\delta_{C_k}((a^T)^{-np+n_kp+k}x)=(\ast).$$
We isolate the terms that have a non-negative power of $a^T$:
$$(\ast)=\left(\prod_{k=0}^{p-1}\prod_{j=1}^{n_k}\hat\delta_{C_k}((a^T)^{(n_k-j)p+k}x)\right)\cdot\prod_{k=0}^{p-1}\prod_{j=n_k+1}^\infty\hat\delta_{C_k}((a^T)^{(n_k-j)p+k}x)=$$
$$\left(\prod_{k=0}^{p-1}\prod_{l=0}^{n_k-1}\hat\delta_{C_k}((a^T)^{lp+k}x)\right)\cdot\prod_{k=0}^{p-1}\prod_{l=1}^\infty\hat\delta_{C_k}((a^T)^{-lp+k}x)$$

Take $0=n_0=\dots=n_{p-1}$ and the computations above show that:
$$\prod_{k=0}^{p-1}\prod_{l=1}^\infty\hat\delta_{C_k}((a^T)^{-lp+k}x)=\hat\mu_{a^p,C_0\oplus
aC_1\oplus\dots\oplus a^{p-1}C_{p-1}}.$$

On the other hand
$$\prod_{k=0}^{p-1}\prod_{l=0}^{n_k-1}\hat\delta_{C_k}((a^T)^{lp+k}x)=\prod_{k=0}^{p-1}\prod_{l=0}^{n_l-1}\hat\delta_{a^{lp+k}C_k}.$$
This product is then the Fourier transform of the convolution of the
measures $\delta_{a^{lp+k}C_k}$. Since the elements in $C_k$ are not
congruent $\mod a\bz^d$, this convolution is $\delta_F$, where $F$ is
as in \eqref{eqf}.

Then the first conclusion follows.

The last statement of our Lemma is now obvious.
\end{proof}

\begin{corollary}\label{cor3_16}
Let $(\tau_b)_{b\in B}$ be an affine IFS in dimension $d=1$, defined by $\tau_b(x)=A^{-1}(x+b)$, $b\in B$. For some $A\in\bz$, $A\geq 1$. Suppose $A=a^p$ for some $a\in\bz$, $p\geq 2$, and 
$$B=a^{n_0p}C\oplus a^{n_1p+1}C\oplus\dots\oplus a^{n_{p-1}p+p-1}C,\mbox{ where }C=\{0,\dots,a-1\},$$
for some non-negative $n_0,\dots,n_{p-1}\in\bz$
Then the invariant measure $\mu_B$ has the form $\mu_B=\mu_{[0,1]}\ast\delta_F$ where $\mu_{[0,1]}$ is the Lebesgue measure on $[0,1]$ and $F$ is a finite set of integers.
\end{corollary}

\begin{proof}
The corollary follows from Lemma \ref{lemdecomp} and the fact that $\mu_{a,C}=\lambda_{[0,1]}$.
\end{proof}

\begin{example}\label{ex3_17}
Consider the affine IFS given by $A=4$, $B=\{0,1,8,9\}$ as in Remark \ref{remlw}. Then $A=2^2$ and 
$B=2^{2\cdot0}\{0,1\}\oplus 2^{2\cdot 1+1}\{0,1\}$. In this case the set $F=\{0\}\oplus 2\{0,1\}=\{0,2\}$.  Therefore we see that $\mu_{A,B}$ is 
$\lambda_{[0,1]}*\delta_{\{0,2\}}$, which is the normalized Lebesgue measure on $[0,1]\cup[2,3]$.

Note also that $B$ has spectrum $\frac{1}{16}\{0,1,8,9\}$.

\end{example}

\begin{theorem}\label{thpro}
Let $\tau_b(x)=A^{-1}(x+b)$, $b\in B$, $x\in\br$, be an affine IFS with $A\in\bz$, $A\geq2$ and $B\in\bz$, $0\in\bz$. Assume the following conditions are satisfied:
\begin{enumerate}
\item[(a)] $A=a^p$ for some $a\in\bz_+$, $p\in\bz_+$. 
\item[(b)] $B=\oplus_{k=0}^{p-1}a^{n_kp+k}C_k$ for some integers $n_0,\dots,n_{p-1}\geq0$ and some subsets $C_0,\dots C_{p-1}$ in $\bz$. 
\item[(c)] For each $k\in\{0,\dots,p-1\}$, the set $C_k$ has spectrum $\frac{1}{a}L_k$ for some subset $L_k$ of $\bz$.
\item[(d)] $\gcd(C_0)=1$.
\end{enumerate}
Then the invariant measure $\mu_B$ is spectral.
\end{theorem}

\begin{proof}
With Lemma \ref{lemdecomp}, $\mu_B$ can be written as the convolution $\mu_{a^p,C_0\oplus aC_1\oplus\dots\oplus a^{p-1}C_{p-1}}\ast\delta_F$, where 
$$F=\oplus_{k=0}^{p-1}\oplus_{l=0}^{n_k-1}a^{lp+k}C_k.$$

With Corollary \ref{cor3_12}, the set $C_0\oplus aC_1\oplus\dots\oplus a^{p-1}C_{p-1}$ has spectrum 
$\frac{1}{a^p}(a^{p-1}L_0\oplus a^{p-2}L_1\oplus\dots\oplus L_{p-1})$. Then, with Theorem \ref{thdj}, the measure 
$ \mu_{a^p,C_0\oplus aC_1\oplus\dots\oplus a^{p-1}C_{p-1}}$ is a spectral measure, with spectrum contained in $\bz$.

From Corollary \ref{cor3_12}, the set $F$ is also spectral. Therefore, the conclusion follows from Corollary \ref{corfinite}.
\end{proof}

\begin{example}\label{ex3_19}
Let us consider the invariant measure $\mu_B$ associated to the affine IFS with $A=4$, $B=\{0,1,4,5\}$.
This example is very similar to the one in Example \ref{ex3_17}, in that $B$ is spectral with spectrum $\frac{1}{8}\{0,1,4,5\}$. However the decomposition of $B$ is $B=\{0,1\}\oplus 2^2\{0,1\}$, so Theorem \ref{thpro} does not apply here. 

\begin{proposition}\label{pr3_24}
The invariant measure $\mu_B$ associated to the affine IFS with $A=4$, $B=\{0,1,4,5\}$ is not a spectral measure. There is however an infinite orthonormal family of exponential functions $e_\lambda$ in $L^2(\mu_B)$.

\end{proposition}

\begin{proof}
Using the decomposition of $B$ and Lemma \ref{lemftinv}, we obtain that $\hat\mu_B=\hat\nu_1\hat\nu_2$, where 
$\nu_1$ is the invariant measure associated to the affine IFS with $A=4$, $B=\{0,1\}$, and $\nu_2$ is the invariant measure associated to the affine IFS with $A=4$, $B=\{0,4\}$. 

We have
$$\hat\nu_1(x)=\prod_{n=1}^\infty\frac{1}{2}(1+e^{2\pi i\frac{x}{4^n}}),\quad
\hat\nu_1(x)=\prod_{n=1}^\infty\frac{1}{2}(1+e^{2\pi i4\frac{x}{4^n}}).$$
The zeros of $\hat\nu_1$ are the points of the form $x=\frac{4^n(2k+1)}2$ with $n\geq 1$, $k\in\bz$. The zeros of $\hat\nu_2$ are the points of the form $x=\frac{4^n(2k+1)}8$ with $n\geq 1$, $k\in\bz$. Thus the zeros of $\hat\nu_1$ are contained in the zeros of $\hat\nu_2$. So the zeros of $\hat\mu_B$ are the same as the zeros of $\hat\nu_2$. 

Suppose $\Lambda$ is a spectrum for $\mu_B$. Then for all $\lambda\neq\lambda'$ in $\Lambda$, $\lambda-\lambda'$ is a zero for $\hat\mu_B$, hence for $\hat\nu_2$. Thus $(e_\lambda)_{\lambda\in\Lambda}$ is an orthonormal family in $L^2(\nu_2)$.

With Lemma \ref{lemspectral} we have
$$1=\sum_{\lambda\in\Lambda}|\hat\mu_B(x+\lambda)|^2=\sum_{\lambda\in\Lambda}|\hat\nu_1(x+\lambda)|^2|\hat\nu_2(x+\lambda)|^2=(\ast).$$
Since $|\hat\nu_1(x)|\leq 1$, with strict inequality for some points $x\in\br$, we get for such points $x$ (exclude the zeros of $\hat\nu_2$):
$$(\ast)<\sum_{\lambda\in\Lambda}|\hat\nu_2(x+\lambda)|^2=\sum_{\lambda\in\Lambda}|\ip{e_{-x}}{e_\lambda}_{L^2(\nu_2)}|^2\leq \|e_{-x}\|^2_{L^2(\nu_2)}=1.$$
We used the orthogonality of $e_\lambda$ for the last inequality, and Lemma \ref{lemspectral} for the last equality.
This contradiction implies that $\mu_B$ cannot be spectral.

To see that there is an infinite orthogonal family of exponentials in $L^2(\mu_B)$, note that $\nu_1$ is a spectral measure because $\{0,1\}$ has spectrum $\frac14\{0,2\}$, so we can use Theorem \ref{thdj}. Since $\hat\mu_B=\hat\nu_1\hat\nu_2$, if $\Lambda$ is a spectrum for $\nu_1$, then $\{e_\lambda\}_{\lambda\in\Lambda}$ is an orthonormal family in $L^2(\mu_B)$.
\end{proof}
\end{example}

\begin{remark}
Earlier work \cite{JoPe98, DJ07a, DJ07b} on Fourier-spectral theory
for IFS-measures $\mu$ on $\br^d$ suggests the following dichotomy: When the
measure $\mu$ is given, then either $L^2(\mu)$ has an orthogonal basis of
Fourier exponentials, or else the set of orthogonal functions $e_\lambda$ in
$L^2(\mu)$ is finite.

   For example \cite{JoPe98}, if $\mu$ is the middle-third Cantor measure with scale
dimension $\log_3(2)$, then $L^2(\mu)$ contains no more than two orthogonal
functions $e_\lambda$. Later work by the present co-authors shows that the
geometry of IFS measures is rather rigid and suggests finiteness of the set
of orthogonal $e_\lambda$ functions in $L^2(\mu)$ unless $\mu$ is in fact a
spectral measure. Hence the dichotomy!

   The conclusion in Proposition \ref{pr3_24} breaks with the dichotomy.
\end{remark}

\begin{theorem}\label{thspec}
Let $A$ be a finite subset of $\bz_+$ with $0\in A$. The following affirmations are equivalent:
\begin{enumerate}
\item $A+[0,1]$ is a spectral set.
\item $A$ is a spectral set.
\end{enumerate}
In this case, any spectrum of $A+[0,1]$ has the form $\bz+\Lambda_A$, where $\Lambda_A$ is a spectrum for $A$. Moreover $A+[0,1]$ has only finitely many spectra that contain $0$.
\end{theorem}

\begin{proof}
(ii)$\Rightarrow$(i) follows from Corollary \ref{corfinite}. This implies also that $\Lambda_A+\bz$ is a spectrum for $A+[0,1]$.

(i)$\Rightarrow$(ii)
Let $\Lambda$ be a spectrum for $[0,1]+A$. We claim that $\Lambda+\bz=\Lambda$.

First, we notice that the Fourier transform of the Lebesgue measure on $[0,1]+A$ is 
\begin{equation}\label{eqhatmua}
\hat\mu_A(t)=\frac{e^{2\pi i t}-1}{2\pi it}(\sum_{a\in A} e^{2\pi ia\cdot t})=:\frac{e^{2\pi it}-1}{2\pi it}p_A(e^{2\pi i t}),\quad(t\in\br\setminus\{0\}) .
\end{equation}

Take $\lambda\in\Lambda$ and $n\in\bz$. Suppose $\lambda+n$ is not in $\Lambda$. Let $\lambda'\in\Lambda$, $\lambda'\neq\lambda$. For the inner product $\ip{\cdot}{\cdot}$ in $L^2(\mu)$, using the fact that $A\subset\bz$ we then have:
$$\ip{e_{\lambda+n}}{e_{\lambda'}}=\hat\mu_A(\lambda+n-\lambda')=\frac{e^{2\pi i(\lambda+n-\lambda')}-1}{2\pi i(\lambda+n-\lambda')}p_A(e^{2\pi i (\lambda-\lambda')})=\frac{\lambda-\lambda'}{\lambda+n-\lambda'}\hat\mu_A(\lambda-\lambda')=0.$$
Since $(e_\lambda)_{\lambda\in\Lambda}$ is an orthonormal basis in $L^2([0,1]+A)$ this shows that $e_{\lambda+n}$ is a scalar multiple of $e_{\lambda}$. This means that $e^{2\pi i n\cdot x}$ must be constant a.e. on $A+[0,1]$. The contradiction implies that $\lambda+n$ must be in $\Lambda$. Thus $\Lambda+\bz\subset\bz$.

Next, let $\Lambda_A:=\Lambda\cap[0,1)$. With the previous facts, we have $\Lambda=\Lambda_A+\bz$. 
For any $\lambda\neq\lambda'\in\Lambda_A$, we have $\hat\mu_A(\lambda-\lambda')=0$. With \eqref{eqhatmua} we obtain 
$p_A(e^{2\pi i(\lambda-\lambda')})=0$. This implies that the rows in the matrix $\frac{1}{\sqrt{N}}(e^{2\pi i a\cdot \lambda})_{\lambda\in\Lambda_A,a\in A}$ are orthogonal. Also this implies that the number of rows in this matrix is smaller or equal to the number of columns, i.e., $\#\Lambda_A\leq\# A$. We have to prove that this matrix is a square matrix. 

For this consider the following: if a function $f\in L^2(A+[0,1])$ is orthogonal to all $(e_\lambda)_{\lambda\in\Lambda}$, then $f=0$. 
We rewrite this orthogonality condition: for all $\lambda\in\Lambda_A$ and $n\in\bz$,
$$0=\ip{e_{\lambda+n}}{f}=\sum_{a\in A}\int_0^1 \cj f(x+a)e^{2\pi i(\lambda+n)\cdot (a+x)}\,dx=
\int_0^1(\sum_{a\in A}\cj f(x+a)e^{2\pi i\lambda\cdot a})e^{2\pi i\lambda\cdot x}e^{2\pi in\cdot x}.$$
Since $(e_n)_{n\in\bz}$ is an orthonormal basis in $L^2([0,1])$, this is equivalent to 
\begin{equation}\label{eqind}
\sum_{a\in A}\cj f(x+a)e^{2\pi i\lambda\cdot a}=0,\quad(\lambda\in\Lambda_A,x\in[0,1]).
\end{equation}
The condition in \eqref{eqind} should be equivalent to $f=0$ on $[0,1]+A$. Thus the vectors
$(e^{2\pi i\lambda\cdot a})_{\lambda\in\Lambda_A}$ must be linearly independent $a\in A$. Then the columns in the matrix $\frac{1}{\sqrt{N}}(e^{2\pi i\lambda\cdot a})_{\lambda\in\Lambda_A,a\in A}$ are linearly independent. This implies that the number of columns is less than the number of rows, i.e., $\#A\leq\#\Lambda_A$. 

So the matrix is square. And the theorem is proved.

Analyzing the proof, we see that any spectrum must have the form $\Lambda_A+\bz$.

Finally, we prove that $A+[0,1]$ has finitely many spectra. Since any spectrum is given by $\Lambda_A+\bz$, with $\Lambda_A$ a spectrum for $A$, it is enough to prove that, $\mod\bz$, there are only finitely many such sets $\Lambda_A$ that contain $0$. Indeed, if $\Lambda_A$ is a spectrum for $A$, then since $A$ is contained in $\bz$, $L:=\Lambda_A\mod\bz$ is a spectrum for $A$ too (because $e^{2\pi iax}=e^{2\pi ia(x\mod\bz)}$ if $a\in\bz$). Note that, from the arguments above, no two elements in $\Lambda_A$ can be congruent $\mod\bz$. If $0\in L$, then for all $l\in L$ one has 
$$\sum_{a\in A}e^{2\pi ial}=0,$$
because of the spectral property of $L$. This implies that $e^{2\pi il}$ is a root of the polynomial $\sum_{a\in A}z^a$. Since this has finitely many roots on the unit circle, and since $L\subset[0,1)$, this shows that there can be only finitely many such sets $L$. This completes the proof.

\end{proof}

\begin{example}
There are sets of the form $A+[0,1]$ that have more than one spectrum. Take $A=\{0,2,4\}$, so $A+[0,1]=[0,1]\cup[2,3]\cup[3,4]$. Then $\{0,\frac13,\frac23\}$ and $\{0,\frac16,\frac26\}$ are spectra for $A$. 
So $\bz+\{0,\frac13,\frac23\}=\frac13\bz$ and $\bz+\{0,\frac16,\frac26\}$ are spectra for $A+[0,1]$.
\end{example}

\begin{theorem}\label{thtile}\cite{New77}
Let $A$ be a finite subset of the integers with $0\in A$. Then $A+[0,1]$ tiles $\br$ if and only if there exists a set of integers $B$, and some $n\in\bn$ such that $A\oplus B=\bz_n$ (as subsets of $\bz_n$). Any tile set $\mathcal T$ is of the form $B\oplus n\bz$, with $B$ as above.
\end{theorem}

\begin{proof}
If $A+[0,1]$ tiles $\br$ with tile set $\mathcal T$, then $A+[0,1]+\mathcal T=\br$. Therefore $A\oplus\mathcal T=\bz$. 
The result in \cite{New77} shows that the tile set $\mathcal T$ has to be periodic, that is $C+n=C$ for some $n\in\bz$. 
Then $\mathcal T$ is a union of congruence classes modulo $n$, so $\mathcal T=B\oplus n\bz$ for some set $B$. Since $A\oplus\mathcal T=\bz$, this implies that $A\oplus B=\bz_n$ modulo $n$.

Conversely, if $A\oplus B=\bz_n$ modulo $n$, then, in $\bz$, $A\oplus B=C$ and $C$ is a subset of $\bz$ congruent to $\{0,\dots,n-1\}$ modulo $n$. Therefore $C\oplus n\bz=\bz$. So $A\oplus (C\oplus n\bz)=\bz$, and $C\oplus n\bz$ is a tile set for $A+[0,1]$.
\end{proof}

\subsection{Unions of intervals as affine IFSs}
 The main results in this section are Theorems \ref{thtile} and \ref{thfinsp} where we identify atomic spectral pairs and rules for their assembly into molecular configurations of new spectral pairs. This is motivated by earlier work \cite{BrJo99}. We use our geometric composition rules in deriving detailed spectral data for the composite spectral-pair systems.

The paper \cite{Lo67} contains some of what we were doing in part of this section. We are including details here nonetheless for the benefit of the reader, and in order to stress what is needed for our purpose.

\begin{theorem}\label{thifs}
Suppose $A$ is a finite set of integers, $0\in A$. Let $\mathcal A:=A+[0,1]$. The following affirmations are equivalent.

\begin{enumerate}
\item There exists an affine IFS of the form $\tau_b(x)=\frac1n(x+b)$ with $b\in B\subset\bz$ and $n\in\bn$, such that 
$A+[0,1]$ is the attractor of the IFS $(\tau_b)_{b\in B}$, with no overlap.
\item There exists a finite set $C\subset\bz$ such that $A\oplus C=\{0,\dots,n-1\}$ for some $n\in\bz$.
\end{enumerate}
Moreover, with the $C$ and $n$ from (ii), $\mathcal A$ tiles $\br$ by $C+n\bz$. In this case, $A+[0,1]$ is a spectral set. 
\end{theorem}

\begin{proof}
(i)$\Rightarrow$(ii). Note that for any $p\geq 1$, we have $$\tau_{b_0}\dots\tau_{b_{p-1}}x=\frac1{n^p}(x+b_0+nb_1+\dots+n^{p-1}b_{p-1}),\quad(x\in\br,b_0,\dots,b_{p-1}\in B).$$
So the IFS $(\tau_{b_{p-1}}\dots\tau_{b_0})_{b_0,\dots,b_{p-1}\in B}$ has the same form. Thus, by picking $p$ large enough, and replacing the IFS $(\tau_b)_{b\in B}$ by the IFS $(\tau_{b_{p-1}}\dots\tau_{b_0})_{b_0,\dots,b_{p-1}\in B}$, which has the same attractor, we can assume that for all $b\in B$, $\mbox{diam}(\tau_b(\mathcal A))<1$. Since the set $\tau_b(\mathcal A)$ has diameter less than 1, it cannot intersect two connected components of $\mathcal A$, because the gap between them is at least $1$

Pick some connected component of $\mathcal A$. Denote it by $[m_0,n_0]$, $m_0,n_0\in\bz$. Since $[m_0,n_0]\subset\mathcal A=\cup_{b\in B}\tau_b(\mathcal A)$, we have that $[m_0,n_0]$ is the union of all sets $\tau_b(\mathcal A)$ which are contained in it. Thus
$$[m_0,n_0]=\bigcup_{b\in B_0}\tau_b(\mathcal A),$$
for some subset $B_0$ of $B$, 
and the union is disjoint because of the non-overlap. Then
$$[nm_0,nn_0]=\bigcup_{b\in B_0}(\mathcal A+b),\mbox{ so }
[0,n(m_0-n_0)]=\bigcup_{b\in B_0}(\mathcal A+b-nn_0).$$
This implies that $A\oplus (B_0-nn_0)=\{0,\dots,n(m_0-n_0)-1\}.$

(ii)$\Rightarrow$(i). The hypothesis implies that 
$\bigcup_{c\in C}(\mathcal A+c)=[0,n]$, disjoint union. Therefore,
$$\mathcal A=\bigcup_{a\in A}([0,1]+a)=\bigcup_{a\in A}\bigcup_{c\in C}\frac{1}{n}(\mathcal A+c+na),$$
and the union is disjoint. This implies (i).

Since $A+C$ is a tile of $\{0,\dots,n-1\}$ the last statement follows. The fact that $A$ is spectral follows from Theorem \ref{thfinsp}.
\end{proof}

 The following result (Theorem \ref{thfinsp}) is closely related to one in
\cite{PW01}, but we include the details here since our techniques are different.
Specifically, we stress the twisted tensor product of Hadamard matrices
\eqref{eq3_5}; a computational feature motivated by fast Fourier transform
algorithms for finite groups.

\begin{theorem}\label{thfinsp}
Let $A$ be a subset of $\bz_+$ such that there exists $B\subset\bz$ and $n\in\bn$ with $A\oplus B=\{0,\dots, n-1\}$. Then $A$ is a spectral set.
\end{theorem}

\begin{remark}
The sets $A\subset\bz_+$ such that $A\oplus B=\{0,\dots,n-1\}$ were completely classified in \cite{Lo67}. The classification is based on the next two Lemmas. We include here the details for the benefit of the reader and to stress what is needed for our purpose. In addition, the proofs will provide a way to construct the spectrum of the set $A$ by means of tensor products of finite Fourier transform matrices. 
\end{remark}

\begin{proof} We will need some Lemmas.

\begin{lemma}\cite{Lo67}\label{lemind}
Suppose $A\oplus B=\{0,\dots,n-1\}$ with $A,B$ finite subsets of $\bz$ and $0\in A\cap B$. Then one, and only one of the following statements is true:
\begin{enumerate}
\item $A=\{0\}$ or $B=\{0\}$.
\item $1\in A$ and there exists a number $d\geq 2$ which divides $n$ and two sets $C,D$ of integers such that:
\begin{enumerate}
\item[(a)] $A=dC+\{0,\dots,d-1\}$.
\item[(b)] $B=dD$.
\item[(c)] $C\oplus D=\{0,\dots,\frac nd-1\}$.
\end{enumerate}
\item $1\in B$ and (a),(b),(c) above hold with the roles of $A$ and $B$ interchanged.
\end{enumerate}

\end{lemma}

\begin{proof}
Suppose $A,B\neq\{0\}$. 
If for some $c\in\{0,\dots,n-1\}$ we have $c=a+b$ with $a\in A$, and $b\in B$, we say that $c=a+b$ is the decomposition of $c$.

Since $A\oplus B$ contains $1$, the element $1$ is exactly in one of $A$ or $B$. Suppose it is in $A$. If it is in $B$ then we interchange $A$ and $B$. Let $d$ be the smallest non-zero element of $B$. We have $d\geq 2$.

Any number $l$ between $1$ and $d-1$ is in $A+B$. Since $d$ is the smallest non-zero element in $B$ it follows that, the decomposition of $l$ is $l=l+0$, so $l$ is in $A$. Therefore $l$ is not in $B$.

Thus $A\cap\{0,\dots,d-1\}=\{0,\dots,d-1\}$ and $B\cap\{0,\dots,d-1\}=\{0\}$.

Inductive hypothesis: assume that for some $p\geq 1$ we have that $A\cap\{0,\dots,pd-1\}=A_pd+\{0,\dots,d-1\}$ and $B\cap\{0,\dots,d-1\}=B_pd$ for some sets $A_p, B_p\subset\bz$.

We claim that the intersections of $A$ and $B$ with $\{0,\dots,(p+1)d-1\}$ has a similar form.

Take $pd$ if $pd\geq n$, we are done, because no number between $pd$ and $pd+d-1$ can be in $A\cup B$. Suppose $pd<n$.
Let $pd=a+b$ be the decomposition of $pd$. 

Case I. $a,b\neq0$. Then, since $b<pd$, using the inductive hypothesis, we have that $b=td$ for some $t\in B_p, t\neq0$ so $a=sd$ for some $s\in A_p$. Then, for $l\in\{0,\dots,d-1\}$ we have that $sd+l$ is in $A$ and $(sd+l)+td$ is the decomposition of $pd+l$. 
This proves that $pd+l$ cannot be in $A\cup B$ (otherwise, the decomposition will be of the form $a+0$ or $0+b$).
Therefore $A\cap\{0,\dots,(p+1)d-1\}=dA_p+\{0,\dots,d-1\}$ and $B\cap\{0,\dots,(p+1)d-1\}=dB_p$.

Case II. $b=0$. Then $pd\in A$. Let $l\in\{1,\dots,d-1\}$. Suppose $pd+l=a+b$ with $b\neq 0$. Then $pd=a+td$ for some $t\in B_p$, $t\neq 0$. Then $A\ni a=(p-t)d+l$. By the induction hypothesis this implies that $(p-t)d\in A$. Then we have the decomposition $pd=(p-t)d+td$, and this contradicts the Case. Thus $b=0$ and $pd+l$ is in $A$. This implies also that $pd+l$ is not in $B$. Therefore 
$A\cap\{0,\dots,(p+1)d-1\}=d(A_p\cup\{p\})+\{0,\dots,d-1\}$ and $B\cap\{0,\dots,(p+1)d-1\}=dB_p$.

Case III. $a=0$. Then $pd\in B$. Let $l\in\{1,\dots,d-1\}$. Then we have the decomposition $pd+l=l+pd$ (since we know $l\in A$). This implies that $pd+l$ cannot be in $A\cup B$ (otherwise we have a decomposition of the form $a+0$ or $0+b$). Therefore $A\cap\{0,\dots,(p+1)d-1\}=dA_p+\{0,\dots,d-1\}$ and $B\cap\{0,\dots,(p+1)d-1\}=d(B_p\cup\{p\})$.

The induction step is proved. Taking $p$ large enough (so that $pd\geq n$) we obtain that
$A=dA_p+\{0,\dots,d-1\}$ and $B=dB_p$ for some sets of integers $A_p,B_p$. Let $a_0=\max A_p$, $b_0=\max B_p$. Then, since $A+B=\{0,\dots, n-1\}$ we must have
$da_0+d-1+db_0=n-1$ so $d(a_0+b_0)=n$. Thus, $d$ is a divisor of $n$ and $a_0+b_0=n/d$.

Also we have $\{0,\dots,n-1\}=A\oplus B=(dA_p+\{0,\dots,d-1\})+dB_p=d(A_p+B_p)+\{0,\dots,d-1\}$. This implies that we must have $A_p\oplus B_p=\{0,\dots,\frac n{d}\}$.

\end{proof}

  When our spectral pairs can be associated with Hadamard matrices, it is natural to ask for the operation on Hadamard matrices which is induced by composition of spectral pair systems under the sum-operation. The next lemma and the final steps in the proof of Theorem \ref{thfinsp} show that the operation on the Hadamard matrices is a twisted tensor product, modeled on the tensor factorizations going into computation of fast Fourier transforms on finite groups, see e.g., \cite{CA06, LvB07}.
  
\begin{lemma}\cite{Lo67}\label{lemops}
If $A\subset\bz_+$ and $A\oplus B=\{0,\dots,n-1\}$, then $A$ can be obtained from a set of the form $\{0,\dots,c-1\}$ after applying several times the following operations:
\begin{enumerate}
\item[(I)] $C\mapsto dC+\{0,\dots,d-1\}$ for some $d\geq 2$;
\item[(II)] $C\mapsto dC$ for some $d\geq 2$.
\end{enumerate}
\end{lemma}
\begin{proof}
To prove Lemma \ref{lemops} we use Lemma \ref{lemind} inductively: $A$ is either $\{0,\dots,n-1\}$ or can be obtained from a set $A_1$ with $A_1\oplus B_1=\{0,\dots,\frac n{d_1}-1\}$ by applying one of the operations (I) or (II). Then the same procedure can be applied to $A_1$. The algorithm stops when $A_k=\{0,\dots, n_k-1\}$ for some $n_k\in\bn$. 
\end{proof}

We can prove now Theorem \ref{thfinsp} using induction and Lemma \ref{lemops}.

If $A=\{0,\dots, n-1\}$ then one can take $L:=\frac1n\{0,\dots,n-1\}$ and all the conditions are satisfied.

Let $C$ be a set in $\bz_+$ such that there exists a set $L$ with $0\in L$, $\#L=\#C=N$ and the matrix 
$\frac{1}{\sqrt{N}}(e^{2\pi ic\cdot l})_{c\in C,l\in L}$ is unitary. We check that the sets obtained by applying the operations (I) and (II) to $C$ has the same properties.

For $dC$ one can take $L'=\frac{1}{d}L$.

For $dC+\{0,\dots,d-1\}$ one can take $L'=\frac{1}{d}(L+\{0,\dots,d-1\})$. 
The corresponding matrix is, with $N=\#C$:
\begin{equation}\label{eq3_5}
\frac{1}{\sqrt{dN}}(e^{2\pi i(dc+k)\cdot\frac1d(l+j)})_{(c,k),(l,j)}=\frac{1}{\sqrt{N}}(e^{2\pi ic\cdot l})_{c\in C,l\in L}\otimes\frac{1}{\sqrt{d}}(e^{2\pi i\frac1d k\cdot l})_{k,l\in\{0,\dots,d-1\}},\end{equation}
so it is unitary. The second matrix in the tensor product is the Fourier transform on the finite group $\bz_d$.

The statement of the theorem follows now by induction and Lemma \ref{lemops}.

For the last statement, suppose $l-l'\in\bz$ for some $l,l'\in L$ with $l\neq l'$. Then $e^{2\pi ia\cdot l}=e^{2\pi ia\cdot l'}$ for all $a\in A$. This shows that the $l$-th and $l'$-th rows in the unitary matrix are equal, and this is a contradiction.

\end{proof}

\subsection{New spectra from old}

          The main result in this section is Theorem \ref{thinf}: For the composite spectral-pair systems in the previous section we prove that a fixed spectral-pair has an infinite set of different spectra.

\begin{lemma}\label{lemnews} Let $\mu_B$ be the invariant measure associated to the IFS 
$\tau_b(x)=A^{-1}(x+b)$, $b\in B$, $x\in\br^d$. Consider the following assumptions on the invariant measure $\mu_B$:
\begin{enumerate}
\item[(a)] The measure $\mu_B$ has spectrum $\Lambda$.
\item[(b)] For every $\lambda\in\Lambda$, $\hat\delta_B(x+\lambda)=\hat\delta_B(x)$ for all $x\in\br^d$. 
\end{enumerate}
Then
\begin{enumerate}
\item  If the finite set $B$ is spectral with spectrum $(A^T)^{-1}L$ for some $L\subset \br^d$,
then $A^T\Lambda\oplus L$ is also a spectrum for $\mu_B$.
\item If there exists a set $L$ in $\br^d$ such that $A^T\Lambda\oplus L=\Lambda$, then $(A^T)^{-1}L$ is a spectrum for $B$.
\end{enumerate}
\end{lemma}

\begin{proof} (a)
With Lemma \ref{lemftinv}, we have for all $x\in\br^d$:
\begin{equation}\label{eqnlam}
\sum_{\lambda\in\Lambda}|\hat\mu_B(A^T(x+\lambda))|^2=\sum_{\lambda\in\Lambda}|\hat\delta_B(x+\lambda)|^2|\hat\mu_B(x+\lambda)|^2=|\hat\delta_B(x)|^2\sum_{\lambda\in\Lambda}|\hat\mu_B(x+\lambda)|^2=|\hat\delta_B(x)|^2.
\end{equation}
We used Lemma \ref{lemspectral} for the last equality.

Since $(A^T)^{-1}L$ is a spectrum for $B$, we have, using Lemma \ref{lemfini} on $\delta_B$:
\begin{equation}\label{eqspdelta}
\sum_{l\in L}|\hat\delta_B(x+(A^T)^{-1}l)|^2=1,\quad(x\in\br^d).
\end{equation}

Using \eqref{eqnlam} and \eqref{eqspdelta},
\begin{equation}\label{eqnlam2}
\sum_{l\in L}\sum_{\lambda\in\Lambda}|\mu_B(A^T(x+(A^T)^{-1}l+\lambda)|^2=\sum_{l\in L}|\hat\delta_B(x+(A^T)^{-1}l)|^2=1.
\end{equation}

Making the substitution $A^Tx=y$, we get
$$\sum_{l\in L,\lambda\in\Lambda}|\hat\mu_B(y+l+A^T\lambda)|^2=1,\quad(y\in\br^d).$$
First, this shows that the writing of an element $a$ as $a=l+A^T\lambda$ with $l\in L$ and $\lambda\in\Lambda$ is unique. Otherwise, take $y=-a$, and on the lefthand side, the sum is $\geq2$, because $\hat\mu_B(0)=1$.
With this, and Lemma \ref{lemspectral}, $L\oplus A^T\Lambda$ is a spectrum.

For (b) we can use the first equality in \eqref{eqnlam2}. Now, the lefthand side of this equality is equal to $1$ since $A^T\Lambda\oplus L=\Lambda$ is a spectrum. Therefore , for all $x\in\br^d$,
$$\sum_{l\in L}|\hat\delta_B(x+(A^T)^{-1}l)|^2=1.$$
With Lemma \ref{lemfini}, this shows that $(A^T)^{-1}L$ is a spectrum for $B$.
\end{proof}

\begin{example}
Let $\mu_4$ be the invariant measure associated to the affine IFS $\tau_{b}(x)=A^{-1}(x+b)$, $b\in B$, with $A=4$ and $B=\{0,2\}$, as in Example \ref{exjp}. We saw that $\Lambda=\{\sum_{k=0}^n4^kl_k\,|\,l_k\in\{0,1\}\}$ is a spectrum for $\mu_4$ \cite{JoPe98}. It is easy to see that for any $q$ odd $L=\frac{1}{4}\{0,q\}$ is a spectrum for $B$. Then, applying Lemma \ref{lemnews} several times one sees that, for any $p\geq0$, the set 
$$\Lambda_{p,q}:=\left\{\sum_{k=0}^n4^kl_k\,|\, l_0,\dots,l_p\in\{0,q\}, l_{p+1},l_{p+2},\dots\in\{0,1\}\right\}$$
is a spectrum for $\mu_4$.
\end{example}

\begin{theorem}\label{thinf}
Let $\mu_B$ be the invariant measure associated to the affine IFS $\tau_b(x)=A^{-1}(x+b)$, $x\in\br$, $b\in B$, where $B$ is a finite set of integers, $0\in B$, and $A\in\bz$, $A\geq 2$. Suppose there exists a set $L$ of integers with $0\in L$ such that $A^{-1}L$ is a spectrum for $B$. If $\#B<A$, then $\mu_B$ has infinitely many spectra. 
\end{theorem}

\begin{proof}
We can assume that $\gcd(B)=1$. If not then let $D:=\gcd(B)$. Let $B':=\frac{1}{D}B$. It is easy to see that $A^{-1}DL$ is a spectrum for $B'$. Also, for a continuous compactly supported function on $\br$,
$$\int f\,d\mu_B=\int f(Dx)\,d\mu_{B'}(x).$$
This implies that $\Lambda$ is a spectrum for $\mu_B$ iff $D\Lambda$ is a spectrum for $\mu_{B'}$ (Corollary \ref{corr}).

Thus we may assume that $\gcd(B)=1$. Also, we can assume that $L$ is contained in $\{0,\dots, A-1\}$. To see this, note that if $L'$ is congruent to $L$ $\mod A$, then $A^{-1}L'$ is spectrum for $B$ too (because $B$ is in $\bz$). Thus, we may replace $L$ by $L\mod A$. We remark that, the fact that $A^{-1}L$ is a spectrum for $B$ implies that no two elements of $L$ are congruent $\mod A$.

First let us analyze the spectrum $\Lambda$ of $\mu_B$ given by Theorem \ref{thdj} as in \cite{DuJo06}. Recall that a $\hat\delta_B$-cycle is a finite set $C:=\{x_0,x_1,\dots x_{p-1}\}$ such that $\tau_{l_0}x_0=x_1,\dots,\tau_{l_{p-2}}x_{p-2}=x_{p-1}$ and $\tau_{l_{p-1}}x_{p-1}=x_0$ for some $l_0,\dots, l_{p-1}\in L$, and $|\hat\delta_B(x_i)|=1$ for $i\in\{0,\dots,p-1\}$.
 (Here $\tau_{l}(x)=A^{-1}(x+l)$ for $x\in\br$, $l\in L$).

As in the proof of Theorem \ref{thdj}, $C$ must be contained in $\bz$. Moreover, since $L\subset\{0,\dots, A-1\}$, we must have that $C\subset[0,1]$. Thus the only possible cycles are $\{0\}$ and $\{1\}$. 

The cycle $\{0\}$ will contribute with 
$$\Lambda(0):=\left\{\sum_{k=0}^nA^kl_k\,|\, l_k\in L, n\geq 0\right\}$$
to the spectrum $\Lambda$.  

If $\{1\}$ is a cycle, then $A-1$ is in $L$, $1=\tau_{A-1}(1)$. It will contribute with 
$$\Lambda(1):=\left\{-A^{n+1}+\sum_{k=0}^nA^{k}l_k\,|\,l_k\in L, n\geq 0\right\}$$
to the spectrum $\Lambda$. 
Thus $\Lambda\subset\Lambda(0)\cup\Lambda(1)$. See \cite{DuJo06} for details.

Note that $\Lambda(1)$ contains only negative numbers. 

Next we will construct a sequence of numbers $n_i\geq0$ in $\bz$ such that $n_{i}-n_{j}\not\in\Lambda$ for all $i>j$.
Since $\#L=\#B<A$, there is an integer $\cj l\in \{0,\dots, A-1\}\setminus L$. Let $n_i=\sum_{k=0}^iA^k\cj l$. If $i>j$, then $n_i-n_j=\sum_{k=j+1}^iA^k\cj l$. This is not in $\Lambda(0)$, because it has base $A$ expansion containing the digit $\cj l$, and it is not in $\Lambda(1)$ since it is positive. Thus, $n_i-n_j$ is not in $\Lambda$.

Pick some fixed element $l_0\neq 0$ in $L$. Define $L_i:=(L\setminus\{l_0\})\cup\{An_i+l_0\}$. Then $L_i$ is a spectrum for $B$ since $L_i\equiv L$ $\mod A$. By Lemma \ref{lemnews}, $A\Lambda\oplus L_i$ is a spectrum for $\mu_B$.

We claim that $A\Lambda\oplus L_i$ are distinct. Suppose $A\Lambda\oplus L_i=A\Lambda\oplus L_j$ for some $i>j$. Then 
$An_i+l_0\in A\Lambda\oplus L_j$. This implies $An_i+l_0=A\lambda+l$ for some $\lambda\in\Lambda$ and $l\in L_j$. But then $l=An_j+l_0$ (since the elements of $L_j$ are distinct $\mod A$). This implies $n_i-n_j\in\Lambda$, a contradiction.

Thus the measure $\mu_B$ has infinitely many spectra.

\end{proof}

\section{Higher dimensions}\label{higher}

       We now return to Fuglede's spectral problem \cite{Fug74} for subsets $\Omega$ in $\br^d$ of finite positive Lebesgue measure. We are concerned with choices of subsets $\Lambda$ such that the associated functions $(e_\lambda)$ indexed by points in $\Lambda$ form an orthogonal basis (ONB) for the Hilbert space $L^2(\Omega)$ with the measure on $\Omega$ being the restriction of $d$-dimensional Lebesgue measure. By analogy to Fourier series we say that the points in $\Lambda$ are Fourier frequencies. When such a choice is possible, we say that the two sets $(\Omega, \Lambda)$ form a spectral pair in $\br^d$.

  In this section we also discuss the variety of possibilities in a spectral pair, when one of the two sets in the pair is fixed. For example, for a spectral pair $(\Omega,\Lambda)$ in $\br^d$, if $\Omega$ is fixed, what are the possibilities for $\Lambda$? And the analogous question when a particular spectrum $\Lambda$ is fixed and given.

For particular cases of sets $\Omega$, the possible spectra $\Lambda$ are known; see \cite{JoPe99} and Theorem \ref{thspec} above for quasi-periodic cases of spectral pairs.
\begin{problem}

Given  $\Omega$, or $\mu$ in one of the families studied in sections 3 and 4 below, write down the structure of the following sets:
$$\{\Lambda\, |\,  (\Omega, \Lambda)\mbox{ is a spectral pair }\},\quad  \{\Lambda\, |\,  (\mu, \Lambda)\mbox{ is a spectral pair }\};\mbox{ and }\{\mathcal T\, |\,  (\Omega, \mathcal T)\mbox{ is a translation pair }\}.$$
(See Definition \ref{defsp}.)
\end{problem}

The answer is known to this when $\Omega$ is the $d$-cube, see \cite{JoPe99, IP98}. For other examples, see also \cite{IKT03}.

Similarly, in case a suitable set $\Lambda$ is specified, the set of measures $\mu$ given by $\{\mu\, |\,  (\mu, \Lambda)\mbox{ is a spectral pair }\}$ would have some interesting structure.

       Because of applications to the study of commuting differential operators, the initial spectral problem considered by Fuglede in \cite{Fug74} concerned the possibilities of spectral pairs in $\br^d$ of the form $(\Omega, \Lambda)$ in $\br^d$ where the set $\Omega$ is assumed {\it connected}.

        A rank-$d$ lattice is a rank-$d$ (discrete) subgroup in $\br^d$. The first observation in \cite{Fug74} was that, if $\Omega$ is a fundamental domain for a lattice $\Gamma$, then $(\Omega, \Lambda)$ is a spectral pair when we take $\Lambda$ to be the lattice dual to $\Gamma$. These spectral pairs are said to be of {\it lattice type}.

      The search for a richer family of spectral pairs, not based on lattices for $d = 1$, leads to the classes of spectral pairs from section \ref{spectral} above. This is of relevance to quasi-periodic structures in solid state physics \cite{BaMo00,BaMo01}. However, note that in our one-dimensional examples, we build sets $\Omega$ as the union of intervals. Or more generally, our spectral pairs have connected components that serve as atoms for composite spectral pairs; see Theorem \ref{thifs} above. Restricting now attention to open sets $\Omega$ in $\br^d$, $d > 1$, we will say that a subset $\Omega$ in $\br^d$ is disconnected if it is the union of connected components, with different components having disjoint closures.

      It is still an open question whether there are any spectral pairs $(\Omega, \Lambda)$ in $\br^2$ with $\Omega$ connected but which are not of lattice type. In this section we will give a procedure for inducing from 1D spectral pairs to 3D spectral pairs in such a way that an induced 3D spectral pair $(\Omega, \Lambda)$ will have $\Omega$ {\it connected} even if the associated 1D set $\Omega$ is disconnected.

\begin{proposition}\label{propprod}
Let $\mu_1$ be a probability measure on $\br^{d_1}$ with spectrum $\Lambda_1$, and suppose for $\mu_1$-a.e. $x_1\in\supp{\mu_1}$, $\mu_{2,x_1}$ is a probability measure on $\br^{d_2}$ with spectrum $\Lambda_2$ (independent of $x_1$). Define the measure $\mu$ on $\br^{d_1+d_2}$ by
$$\int_{\br^{d_1+d_2}}f(t_1,t_2)\,d\mu(t_1,t_2)=\int_{\br^{d_1}}\int_{\br^{d_2}}f(t_1,t_2)\,d\mu_{2,t_1}(t_2)\,d\mu_1(t_1),\quad(f\in C_c(\br^{d_1+d_2})).$$
Then $\mu$ is a spectral measure with spectrum $\Lambda_1\times\Lambda_2$.
\end{proposition}

\begin{proof}
We have for all $x_1\in\br^{d_1},x_2\in\br^{d_2}$,
\begin{equation}\label{eqmu1}
\hat\mu(x_1,x_2)=\int_{\br^{d_1}}e^{2\pi x_1\cdot t_1}\int_{\br^{d_2}}e^{2\pi ix_2\cdot t_2}\,d\mu_{2,t_1}(t_2)\,d\mu_1(t_1)=\int_{\br^{d_1}}e^{2\pi i t_1\cdot x_1}\hat\mu_{2,t_1}(x_2)\,d\mu_1(t_1).
\end{equation}

Then 
$$\sum_{\lambda_1\in\Lambda_1,\lambda_2\in\Lambda_2}|\hat\mu(x_1-\lambda_1,x_2-\lambda_2)|^2=\sum_{\lambda_2,\lambda_1}\left|\int_{\br_{d_1}}e^{2\pi i(x_1-\lambda_1)\cdot t_1}\hat\mu_{2,t_1}(x_2-\lambda_2)\,d\mu_1(t_1)\right|^2=(\ast)$$
But, since $(e_{\lambda_1})_{\lambda_1\in\Lambda_1}$ is an ONB for $L^2(\mu_1)$, the Parseval identity implies that
$$\sum_{\lambda_1}\left|\int_{\br^{d_1}}e^{2\pi i(x_1-\lambda_1)\cdot t_1}\hat\mu_{2,t_1}(x_2-\lambda_2)\,d\mu_1(t_1)\right|^2=\int_{\br^{d_1}}\left|e^{2\pi ix_1\cdot t_1}\hat\mu_{2,t_1}(x_2-\lambda_2)\right|^2\,d\mu_1(t_1)=$$$$\int_{\br^{d_1}}|\hat\mu_{2,t_1}(x_2-\lambda_2)|^2\,d\mu_1(t_1),$$
for all $\lambda_2\in\Lambda_2$.
Therefore, with Lemma \ref{lemspectral},
$$(\ast)=\sum_{\lambda_2}\int_{\br^{d_1}}|\hat\mu_{2,t_1}(x_2-\lambda_2)|^2\,d\mu_1(t_1)=\int_{\br^{d_1}}\sum_{\lambda_2}|\hat\mu_{2,t_1}(x_2-\lambda_2)|^2\,d\mu_1(t_1)=1.$$
Then the same Lemma implies that $\mu$ has spectrum $\Lambda_1\times\Lambda_2$.

\end{proof}

\begin{corollary}\label{corslice}
Let $C$ be a bounded measurable set in $\br^{d_1+d_2}$. Let $A$ be the projection of $C$ onto $\br^{d_1}$. Assume that $A$ has spectrum $\Lambda_1$. Suppose for Lebesgue a.e. $x_1$ in $A$, the slice $C_{x_1}:=\{x_2\in\br^{d_2}\,|\, (x_1,x_2)\in C\}$ is spectral with spectrum $\Lambda_2$ (independent of $x_1$) and $C_{x_1}$ has Lebesgue measure $c$ (independent of $x_1$). Then the set $C$ has spectrum $\Lambda_1\times\Lambda_2$. 
\end{corollary}

\begin{proof}
The corollary follows from Proposition \ref{propprod} by taking $\mu_1$ to be the normalized Lebesgue measure on $A$, and $\mu_{2,x_1}$ to be the normalized Lebesgue measure on $C_{x_1}$. Then $\mu$ will be the normalized Lebesgue measure on $C$.
\end{proof}

\begin{example}\label{exst}
Let $I^3$ be the unit cube. Let $S$ be the set in Figure \ref{Figstairs}, defined by
$$S:=I^3+\left\{(0,0,0),(1,0,\frac13),(2,0,\frac23),(2,1,\frac33),(1,1,\frac43),(0,1,\frac53),\right.$$$$\left.(0,0,\frac63),(1,0,\frac73),(2,0,\frac83),(2,1,\frac93),(1,1,\frac{10}3),(0,1,\frac{11}3)\right\}.$$

We label these cubes by $(i)$ for $i=1\dots 12$ according to the position in the list (or according to the $z$-coordinate: cube $(i)$ is at height $\frac{i-1}3$). 
\begin{figure}
\setlength{\unitlength}{0.4\textwidth}
\begin{picture}(8.5,2)(-0.9,0)                      
\put(0,0){\includegraphics[bb=88 4 376 148,width=\unitlength]{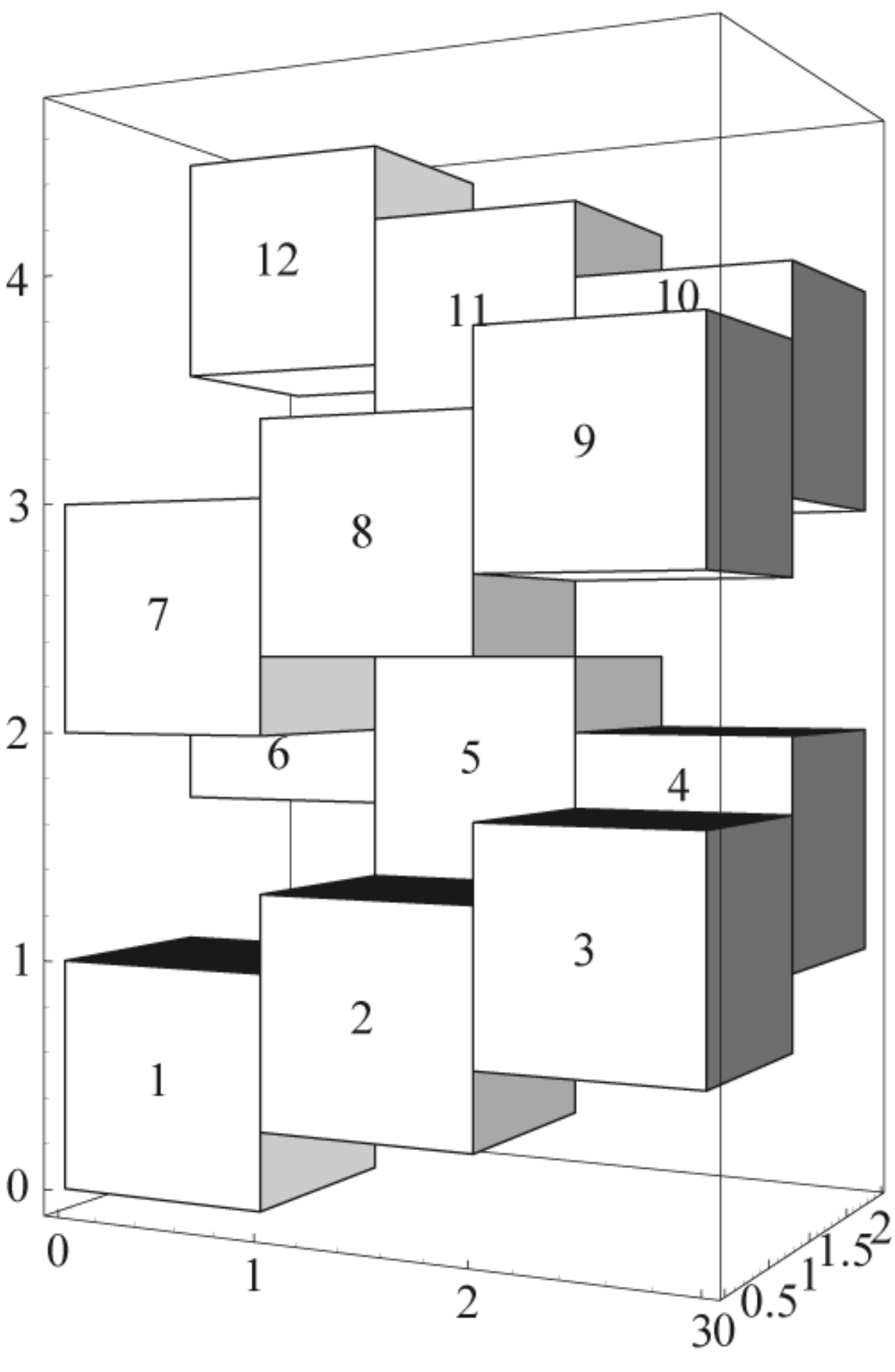}}
\end{picture}
\caption{The staircase}
\label{Figstairs}
\end{figure}

\begin{proposition}
The staircase set $S$ has the following properties:
\begin{enumerate}
\item $S$ is connected.
\item $S$ has spectrum $\frac13\bz\times\frac12\bz\times(\bz+\{0,\frac14\})$.
\item $S$ tiles $\br^3$ by $3\bz\times 2\bz\times (4\bz+\{0,1\})$.
\item $S$ does not tile $\br^3$ by any lattice.
\end{enumerate}
\end{proposition}

\begin{proof}
(i) and (iii) are obvious. For (ii) we use Corollary \ref{corslice}. Note that each $y$-slice $S_y$ is of the form (depending on $y\in[0,1]$ or $y\in[1,2]$):
$$S_y=I^2+\left\{(0,0),(1,\frac13),(2,\frac23),(0,\frac63),(1,\frac73),(2,\frac83)\right\},\mbox{ or }
S_y=I^2+\left\{(2,\frac33),(1,\frac43),(0,\frac53),(2,\frac93),(1,\frac{10}3),(0,\frac{11}3)\right\}.$$

Each $S_y$ has spectrum $\frac13\bz\times(\bz+\{0,\frac14\})$. To see this, we use Corollary \ref{corslice} again, and note that each $x$-slice of $S_y$ is a translation of $[0,1]\cup[2,3]$. From Example \ref{ex3_6}, we know that this has spectrum $\bz+\{0,\frac14\}$. Since $[0,3]$ has spectrum $\frac13\bz$, we obtain with Corollary \ref{corslice} that $S_y$ has spectrum $\frac13\bz\times(\bz+\{0,\frac14\})$. Since $[0,2]$ has spectrum $\frac12\bz$, we obtain (ii).

(iv) Let us assume that $S$ tiles $\br^3$ by some lattice $\Gamma$. Consider the point $(1.5,0.5,1.8)$ located between the cubes $(2)$ and $(8)$. It must belong to one of the translations $S+\gamma$. Thus it belongs to one of the translations of the cubes $(i)+\gamma$, $i\in\{1,\dots,12\}$. Then the cube $(i)+\gamma$ must fit perfectly between the cubes $(2)$ and $(8)$, otherwise there is some space left between the cube $(i)+\gamma$ and $(5)$, and this space has to be covered by another translation of $S$. This would be impossible since the space left between $(5)$ and $(i)+\gamma$ is too small. So the cube $(i)+\gamma$ fits perfectly between the cubes $(2)$,$(8)$ and $(5)$. 

 The cube $(i)$ cannot be one of the cubes $(4),(5),(6),(10),(11),(12)$. Suppose $(i)$ is the cube $(4)$. Then $(5)+\gamma$ intersects the cube $(7)$ which contradicts the tiling property. Suppose $(i)$ is the cube $(12)$. Then $(11)+\gamma$ intersects the cube $(3)$. All cases can be treated in this way. 
 
 Thus the cube $(i)$ is one of the cubes $(1),(2),(3),(7),(8),(9)$.
 Suppose it is the cube $(1)$. We know now that $(1)+\gamma$ fits perfectly between $(2),(8)$ and $(5)$. Thus $\gamma=(1,0,\frac43)$. But the cube $(5)+(1,0,\frac43)=I^3+(2,1,\frac83)$ intersects the cube $(10)=I^3+(2,1,\frac93)$.
  If $(i)$ is the cube $(2)$, then $\gamma=(0,0,1)$. Since $\Gamma$ is a lattice $(0,0,2)$ is also in $\Gamma$
  but $(2)+(0,0,2)$ intersects $S$ in cube $(8)$, so we cannot have tiling.
  
  If $(i)$ is the cube $(3)$, then $\gamma=(-1,0,\frac23)$. Then $2\gamma$ is in $\Gamma$ too. But $(3)+2\gamma$ intersects cube $(7)$, so we cannot have tiling. 
  
 The cases $(7),(8),(9)$ can be treated similarly, and we reach the desired contradiction.

\end{proof}
\end{example}

\begin{example}\label{ex4_5}
Let $C=I^2$ be the standard unit cube in $\br^2$. Divide $C$ along the main diagonal, resulting in two triangles $U$ and $V$, $U$ over the diagonal, and $V$ under. Let $p\in\bz\setminus\{0\}$, and set 
$$\Omega(p):=U\cup ((p,0)+V).$$
Clearly then $\Omega(p)$ is disconnected. Moreover $(\Omega(p),\bz^2)$ is a spectral pair, but neither of the two (separated) components $U$ or $(p,0)+V$ is.
\begin{proof}
Since $(C,\bz^2)$ is trivially a spectral pair, and $\Omega(p)$ is congruent modulo $\bz^2$ to $C$ it follows that $(\Omega(p),\bz^2)$ is one too. But it is clear that the triangles $U$ (or $V$) cannot have any subset $\mathcal L\subset\br^2$ such that $(U,\mathcal L)$ is a spectral pair (see \cite{Fug74}).
\end{proof}

\begin{remark}
The paper \cite{Fug74} contains two important examples of connected open sets $\Omega$ in $\br^2$ for which there is no choice of $\Lambda$ turning the pair $(\Omega, \Lambda)$ into a spectral pair: $\Omega$ the open disk $D$, and $\Omega$ the open triangle $T$. Fuglede showed that for $D$ there cannot be an infinite set of $\lambda$'s which are mutually orthogonal; while for $T$, there are infinite orthogonal families in $L^2(T)$ but none of them are total; hence no orthogonal basis.
\end{remark}

\begin{problem}
Are the only bounded open connected sets in $\br^2$ of finite positive measure which tile $\br^2$ by translations (or are spectral sets) the fundamental domains for rank-2 lattices? In other words, if a bounded open connected subset of $\br^2$ tiles $\br^2$ by translations (or is a spectral set), then does it tile $\br^2$ by some lattice?
\end{problem}

Example \ref{exst} shows that there are examples in $\br^3$ of translation sets
(spectral sets), open and connected of finite positive measure which are not
fundamental domains.

\begin{acknowledgements}
Some of the ideas in the paper grew out of discussions, over the years with the following colleagues, Bent Fuglede, Jeff Lagarias, Gabriel Picioroaga, Steen Pedersen, Michael Reid, and Yang Wang. Especially the examples in section 4 owe much to what we learned from Bent Fuglede and Steen Pedersen. The authors thank the referee for insightful suggestions regarding the presentation, and for adding to the list of references, especially \cite{IKT01} and \cite{Lo67}.
\end{acknowledgements}

\end{example}

\bibliographystyle{alpha}
\bibliography{sifs}

\end{document}